\documentclass[a4paper, 12pt]{article}
\usepackage{amsmath,amsthm,amssymb,graphicx,psfrag}

\def\dom{{\rm dom}}

\def\ol{\bar}

\def\a{\alpha}

\def\rn{\mathbb{R}^n}

\def\re{\mathbb{R}}

\def\dist{{\rm dist}}

\def\sg{\partial}
\def\dom{{\rm dom}}

\newtheorem{theorem}{Theorem}

\newtheorem{assumption}{Assumption}

\newtheorem{lemma}{Lemma}

\newtheorem{proposition}{Proposition}

\begin{document}
\title{Distributed Multi-Agent Optimization with State-Dependent Communication\footnote{ This research was partially
supported by the National Science Foundation under Career grant
DMI-0545910, the DARPA ITMANET program, and the AFOSR MURI
R6756-G2.}}
\author{Ilan
Lobel\footnote{Microsoft Research New England Lab and Stern
School of Business, New York University, ilobel@stern.nyu.edu},
Asuman Ozdaglar\footnote{Laboratory for Information and
Decision Systems, Department of Electrical Engineering and
Computer Science, Massachusetts Institute of Technology,
asuman@mit.edu} and Diego Feijer\footnote{Laboratory for
Information and Decision Systems, Department of Electrical
Engineering and Computer Science, Massachusetts Institute of
Technology, feijer@mit.edu}}

%\markright{LIDS Report ???}

\maketitle

{\it This paper is dedicated to the memory of Paul Tseng, a great
researcher and friend.}

\thispagestyle{empty}

%\vspace{-0.2in}

\begin{abstract}
We study distributed algorithms for solving global optimization
problems in which the objective function is the sum of local
objective functions of agents and the constraint set is given by the
intersection of local constraint sets of agents. We assume that each
agent knows only his own local objective function and constraint
set, and exchanges information with the other agents over a randomly
varying network topology to update his information state. We assume
a {\it state-dependent communication model} over this topology:
communication is Markovian with respect to the states of the agents
and the probability with which the links are available depends on
the states of the agents.

In this paper, we study a {\it projected multi-agent subgradient
algorithm} under state-dependent communication. The algorithm
involves each agent performing a local averaging to combine his
estimate with the other agents' estimates, taking a subgradient step
along his local objective function, and projecting the estimates on
his local constraint set. The state-dependence of the communication
introduces significant challenges and couples the study of
information exchange with the analysis of subgradient steps and
projection errors. We first show that the multi-agent subgradient
algorithm when used with a constant stepsize may result in the agent
estimates to diverge with probability one. Under some assumptions on
the stepsize sequence, we provide convergence rate bounds on a
``disagreement metric" between the agent estimates. Our bounds are
time-nonhomogeneous in the sense that they depend on the initial
starting time. Despite this, we show that agent estimates reach an
almost sure consensus and converge to the same optimal solution of
the global optimization problem with probability one under different
assumptions on the local constraint sets and the stepsize sequence.

\end{abstract}

%\textbf{Keywords:} Constrained multi-agent optimization,
%state-dependent communication.

\newpage

\setcounter{page}{1}

\section{Introduction}

Due to computation, communication, and energy constraints,
several control and sensing tasks are currently performed
collectively by a large network of autonomous agents.
Applications are vast including a set of sensors collecting and
processing information about a time-varying spatial field
(e.g., to monitor temperature levels or chemical
concentrations), a collection of mobile robots performing
dynamic tasks spread over a region, mobile relays providing
wireless communication services, and a set of humans
aggregating information and forming beliefs about social issues
over a network. These problems motivated a large literature
focusing on design of optimization, control, and learning
methods that can operate using local information and are robust
to dynamic changes in the network topology. The standard
approach in this literature involves considering
``consensus-based" schemes, in which agents exchange their
local estimates (or states) with their neighbors with the goal
of aggregating information over an {\it exogenous} (fixed or
time-varying) network topology. In many of the applications,
however, the relevant network topology is configured {\it
endogenously as a function of the agent states}, for example,
the communication network varies as the location of mobile
robots changes in response to the objective they are trying to
achieve. A related set of problems arises when the current
information of decentralized agents influences their potential
communication pattern, which is relevant in the context of
sensing applications and in social settings where disagreement
between the agents would put constraints on the amount of
communication among them.

In this paper, we propose a general framework for design and
analysis of distributed multi-agent optimization algorithms with
state dependent communication. Our model involves a network of $m$
agents, each endowed with a local objective function
$f_i:\mathbb{R}^n\to \mathbb{R}$ and a local constraint
$X_i\subseteq \mathbb{R}^n$ that are private information, i.e., each
agent only knows its own objective and constraint. The goal is to
design distributed algorithms for solving a global constrained
optimization problem for optimizing an objective function, which is
the sum of the local agent objective functions, subject to a
constraint set given by the intersection of the local constraint
sets of the agents. These algorithms involve each agent maintaining
an estimate (or state) about the optimal solution of the global
optimization problem and update this estimate based on local
information and processing, and information obtained from the other
agents.

We assume that agents communicate over a network with randomly
varying topology. Our random network topology model has two novel
features: First, we assume that the communication at each time
instant $k$, (represented by a {\it communication matrix} $A(k)$
with positive entries denoting the availability of the links between
agents) is Markovian on the states of the agents. This captures the
time correlation of communication patterns among the
agents.\footnote{Note that our model can easily be extended to model
Markovian dependence on other stochastic processes, such as channel
states, to capture time correlation due to global network effects.
We do not do so here for notational simplicity.} The second, more
significant feature of our model is that the probability of
communication between any two agents at any time is a function of
the agents' states, i.e., the closer the states of the two agents,
the more likely they are to communicate. As outlined above, this
feature is essential in problems where the state represents the
position of the agents in sensing and coordination applications or
the beliefs of agents in social settings.

For this problem, we study a {\it projected multi-agent subgradient
algorithm}, which involves each agent performing a local averaging
to combine his estimate with the other agents' estimates he has
access to, taking a subgradient step along his local objective
function, and projecting the estimates on his local constraint set.
We represent these iterations as stochastic linear time-varying
update rules that involve the agent estimates, subgradients and
projection errors explicitly. With this representation, the
evolution of the estimates can be written in terms of stochastic
transition matrices $\Phi(k,s)$ for $k\ge s\ge 0$, which are
products of communication matrices $A(t)$ over a window from time
$s$ to time $k$. The transition matrices $\Phi(k,s)$ represent
aggregation of information over the network as a result of local
exchanges among the agents, i.e., in the long run, it is desirable
for the transition matrices to converge to a uniform distribution,
hence aligning the estimates of the agents with uniform weights
given to each (ensuring that information of each agent affects the
resulting estimate uniformly). As a result, the analysis of our
algorithm involves studying convergence properties of transition
matrices, understanding the limiting behavior of projection errors,
and finally studying the algorithm as an ``approximate subgradient
algorithm" with bounds on errors due to averaging and projections.

In view of the dependence of information exchange on the agent
estimates, it is not possible to decouple the effect of stepsizes
and subgradients from the convergence of the transition matrices. We
illustrate this point by first presenting an example in which the
projected multi-agent subgradient algorithm is used with a constant
stepsize $\alpha(k)=\alpha$ for all $k\ge 0$. We show that in this
case, agent estimates and the corresponding global objective
function values may diverge with probability one for any constant
value of the stepsize. This is in contrast to the analysis of
multi-agent algorithms over exogenously varying network topologies
where it is possible to provide error bounds on the difference
between the limiting objective function values of agent estimates
and the optimal value as a function of the constant stepsize $\a$
(see \cite{Ilanoptim}).

We next adopt an assumption on the stepsize sequence\
$\{\a(k)\}$ (see Assumption \ref{ass. small steps}), which
ensures that $\a(k)$ decreases to zero sufficiently fast, while
satisfying $\sum_{k=0}^\infty \a(k)=\infty$ and
$\sum_{k=0}^\infty \alpha^2(k)<\infty$ conditions. Under this
assumption, we provide a bound on the expected value of the
disagreement metric, defined as the difference $\max_{i,j}
|[\Phi(k,s)]_{ij}-{1\over m}|$. Our analysis is novel  and
involves constructing and bounding (uniformly) the probability
of a hierarchy of events, the length of which is specifically
tailored to grow faster than the stepsize sequence, to ensure
propagation of information across the network before the states
drift away too much from each other. In contrast to exogenous
communication models, our bound is {\it time-nonhomogeneous},
i.e., it depends on the initial starting time $s$ as well as
the time difference $(k-s)$. We also consider the case where we
have the assumption that the agent constraint sets $X_i$'s are
compact, in which case we can provide a bound on the
disagreement metric without any assumptions on the stepsize
sequence.

Our next set of results study the convergence behavior of agent
estimates under different conditions on the constraint sets and
stepsize sequences. We first study the case when the local
constraint sets of agents are the same, i.e., for all $i$, $X_i=X$
for some nonempty closed convex set. In this case, using the
time-nonhomogeneous contraction provided on the disagreement metric,
we show that agent estimates reach almost sure consensus under the
assumption that stepsize sequence $\{\a(k)\}$ converges to 0
sufficiently fast (as stated in Assumption \ref{ass. small steps}).
Moreover, we show that under the additional assumption
$\sum_{k=0}^\infty \a(k)=\infty$, the estimates converge to the same
optimal point of the global optimization problem with probability
one. We then consider the case when the constraint sets of the
agents $X_i$ are different convex compact sets and present
convergence results both in terms of almost sure consensus of agent
estimates and almost sure convergence of the agent estimates to an
optimal solution under weaker assumptions on the stepsize sequence.

Our paper contributes to the growing literature on multi-agent
optimization, control, and learning in large-scale networked
systems. Most work in this area builds on the seminal work by
Tsitsiklis \cite{johnthes} and Bertsekas and Tsitsiklis
\cite{distbook} (see also Tsitsiklis {\it et al.}\ \cite{distasyn}),
which developed a general framework for parallel and distributed
computation among different processors. Our work is related to
different strands of literature in this area.

One strand focuses on reaching consensus on a particular scalar
value or computing exact averages of the initial values of the
agents, as natural models of cooperative behavior in
networked-systems (for deterministic models, see \cite{vicsek},
\cite{ali},\cite{reza}, \cite{spielman}, \cite{alexCDC}, and
\cite{alexlong}; for randomized models, where the randomness may be
due to the choice of the randomized communication protocol or due to
the unpredictability in the environment that the information
exchange takes place, see \cite{boyd}, \cite{mesbahi}, \cite{wu},
\cite{alireza-one}, \cite{alireza-two}, and \cite{fagnani}) Another
recent literature studies optimization of more general objective
functions using subgradient algorithms and consensus-type mechanisms
(see \cite{nedic-ozdaglar}, \cite{quantization},
\cite{constconsoptim}, \cite{Ilanoptim}, \cite{baras},
\cite{raminc}, \cite{martinez}). Of particular relevance to our work
are the papers \cite{Ilanoptim} and \cite{constconsoptim}. In
\cite{Ilanoptim}, the authors studied a multi-agent unconstrained
optimization algorithm over a random network topology which varies
independently over time and established convergence results for
diminishing and constant stepsize rules. The paper
\cite{constconsoptim} considered multi-agent optimization algorithms
under deterministic assumptions on the network topology and with
constraints on agent estimates. It provided a convergence analysis
for the case when the agent constraint sets are the same. A related,
but somewhat distinct literature, uses consensus-type schemes to
model opinion dynamics over social networks (see \cite{golub},
\cite{golub-two}, \cite{misinformation}, \cite{krause},
\cite{opinion-dyn}). Among these papers, most related to our work
are \cite{krause} and \cite{opinion-dyn}, which studied dynamics
with opinion-dependent communication, but without any optimization
objective.

The rest of the paper is organized as follows: in Section
\ref{sec:model}, we present the optimization problem, the
projected subgradient algorithm and the communication model. We
also show a counterexample that demonstrates that there are
problem instances where this algorithm, with a constant
stepsize, does not solve the desired problem. In Section
\ref{sec:stochastic}, we introduce and bound the disagreement
metric $\rho$, which determines the spread of information in
the network. In Section \ref{sec:optim}, we build on the
earlier bounds to show the convergence of the projected
subgradient methods. Section \ref{sec:conclusions} concludes.

\vskip 2pc

\noindent {\bf Notation and Basic Relations:} \vskip .5pc

A vector is viewed as a column vector, unless clearly stated
otherwise. We denote by $x_i$ or $[x]_i$ the $i$-th component of a
vector $x$. When $x_i\ge 0$ for all components $i$ of a vector $x$,
we write $x\ge 0$. For a matrix $A$, we write $A_{ij}$ or $[A]_{ij}$
to denote the matrix entry in the $i$-th row and $j$-th column. We
denote the nonnegative orthant by $\mathbb{R}^n_+$, i.e.,
$\mathbb{R}^n_+ = \{x\in \mathbb{R}^n\mid x\ge 0\}$. We write $x'$
to denote the transpose of a vector $x$. The scalar product of two
vectors $x,y\in\re^n$ is denoted by $x'y$. We use $\|x\|$ to denote
the standard Euclidean norm, $\|x\|=\sqrt{x'x}$.

A vector $a\in\re^m$ is said to be a {\it stochastic vector}
when its components $a_i$, $i=1,\ldots,m$, are nonnegative and
their sum is equal to 1, i.e., $\sum_{i=1}^m a_i =1$. A square
$m\times m$ matrix $A$ is said to be a {\it stochastic matrix}
when each row of $A$ is a stochastic vector. A square $m\times
m$ matrix $A$ is said to be a {\it doubly stochastic} matrix
when both $A$ and $A'$ are stochastic matrices.

For a function $F:\re^n\to(-\infty,\infty]$, we denote the
domain of $F$ by $\dom(F)$, where
\[\dom(F)=\{x\in\re^n \mid F(x)<\infty\}.\]
We use the notion of a subgradient of a {\it convex} function
$F(x)$ at a given vector $\ol x\in \dom(F)$. We say that
$s_F(\ol x)\in\re^n$ {\it is a subgradient of the function $F$
at $\ol x\in\dom(F)$} when the following relation holds:
\begin{equation}
F(\ol x) + s_F(\ol x)'(x-\ol x)\le F(x)\qquad \hbox{for all
}x\in\dom(F). \label{sgdconvdef}
\end{equation}
The set of all subgradients of $F$ at $\ol x$ is denoted by $\sg
F(\ol x)$ (see \cite{ourbook}).

In our development, the properties of the projection operation on a
closed convex set play an important role. We write $dist(\ol x, X)$
to denote the standard Euclidean distance of a vector $\ol x$ from a
set $X$, i.e.,
$$\dist(\ol x,X) =\inf_{x\in X}\|\ol x - x\|.$$
Let $X$ be a nonempty closed convex set in $\re^n$. We use
$P_X[\ol x]$ to denote the projection of a vector $\ol x$ on
set $X$, i.e.,
\[P_X[\ol x]=\arg\min_{x\in X}\|\ol x-x\|.\]
We will use the standard non-expansiveness property of projection,
i.e.,
\begin{equation}\label{nonexpan}
\|P_X[x] - P_X[y]\|\le \|x-y\|\qquad\hbox{for any }x \hbox{ and } y.
\end{equation}
We will also use the following relation between the projection error
vector and the feasible directions of the convex set $X$: for any
$x\in \rn$,
\begin{equation}\|P_X[x]-y\|^2\le \|x-y\|^2 - \|P_X[x]-x\|^2 \qquad\hbox{for
all $y\in X$}.\label{projerror-feasdir}\end{equation}

\section{The Model}\label{sec:model}

\subsection{Optimization Model}

We consider a network that consists of a set of nodes (or
agents) $\mathcal{M}=\{1,\dots,m\}$. We assume that each agent
$i$ is endowed with a local objective (cost) function $f_i$ and
a local constraint function $X_i$ and this information is
distributed among the agents, i.e., each agent knows only his
own cost and constraint component. Our objective is to develop
distributed algorithms that can be used by these agents to
cooperatively  solve the following constrained optimization
problem:
\begin{eqnarray}
& \text{minimize}\quad\; \sum_{i=1}^m f_i(x) \label{optim-prob}\\
& \text{subject to}\quad x\in\cap_{i=1}^m X_i,\nonumber
\end{eqnarray}
where each $f_i:\mathbb{R}^n\rightarrow\mathbb{R}$ is a convex
(not necessarily differentiable) function, and each $X_i
\subseteq \rn$ is a closed convex set. We denote the
intersection set by $X=\cap_{i=1}^m X_i$ and assume that it is
nonempty throughout the paper. Let $f$ denote the global
objective, that is, $f(x) = \sum_{i=1}^m f_i(x)$, and $f^*$
denote the optimal value of problem (\ref{optim-prob}), which
we assume to be finite. We also use $X^*=\{x\in X:f(x)=f^*\}$
to denote the set of optimal solutions and assume throughout
that it is nonempty.

We study a distributed multi-agent subgradient method, in which
each agent $i$ maintains an {\it estimate} of the optimal
solution of problem (\ref{optim-prob}) (which we also refer to
as the {\it state of agent $i$}), and updates it based on his
local information and information exchange with other
neighboring agents. Every agent $i$ starts with some initial
estimate $x_i(0)\in X_i$. At each time $k$, agent $i$ updates
its estimate according to the following:

\begin{equation}
x_i(k+1) = P_{X_i}\left[\sum_{j=1}^m a_{ij}(k)x_j(k) -
\alpha(k)d_i(k)\right], \label{eq.update_rule}
\end{equation}
where $P_{X_i}$ denotes the projection on agent $i$ constraint set
$X_i$, the vector $[a_{ij}(k)]_{j\in\mathcal{M}}$ is a vector of
weights for agent $i$, the scalar $\alpha(k)>0$ is the stepsize at
time $k$, and the vector $d_i(k)$ is a subgradient of agent $i$
objective function $f_i(x)$ at his estimate $v_i(k) = \sum_{j=1}^m
a_{ij}(k)x_j(k)$. Hence, in order to generate a new estimate, each
agent combines the most recent information received from other
agents with a step along the subgradient of its own objective
function, and projects the resulting vector on its constraint set to
maintain feasibility. We refer to this algorithm as the {\it
projected multi-agent subgradient algorithm}.\footnote{See also
\cite{constconsoptim} where this algorithm is studied under
deterministic assumptions on the information exchange model and the
special case $X_i=X$ for all $i$.} Note that when the objective
functions $f_i$ are identically zero and the constraint sets
$X_i=\mathbb{R}^n$ for all $i\in\mathcal{M}$, then the update rule
(\ref{eq.update_rule}) reduces to the classical averaging algorithm
for {\it consensus} or {\it agreement} problems, as studied in
\cite{multiagent} and \cite{ali}.

In the analysis of this algorithm, it is convenient to separate
the effects of different operations used in generating the new
estimate in the update rule (\ref{eq.update_rule}). In
particular, we rewrite the relation in Eq.\
(\ref{eq.update_rule}) equivalently as follows:
\begin{eqnarray}
v_i(k) &=& \sum_{j=1}^m a_{ij}(k)x_j(k),\label{convex-comb}\\
x_i(k+1) &=& v_i(k) -\alpha(k)d_i(k) +e_i(k),\label{subgradient-step}\\
e_i(k) &=& P_{X_i}[v_i(k) -\alpha(k)d_i(k)] - \Big(v_i(k)
-\alpha(k)d_i(k)\Big).\label{proj-error}
\end{eqnarray}
This decomposition allows us to generate the new estimate using a
{\it linear update rule} in terms of the other agents' estimates,
the subgradient step, and the projection error $e_i$. Hence, the
nonlinear effects of the projection operation is represented by the
projection error vector $e_i$, which can be viewed as a perturbation
of the subgradient step of the algorithm. In the sequel, we will
show that under some assumptions on the agent weight vectors and the
subgradients, we can provide upper bounds on the projection errors
as a function of the stepsize sequence, which enables us to study
the update rule (\ref{eq.update_rule}) as an approximate subgradient
method.

We adopt the following standard assumption on the subgradients
of the local objective functions $f_i$.

\begin{assumption}(Bounded Subgradients)
The subgradients of each of the $f_i$ are uniformly bounded,
i.e., there exists a scalar $L>0$ such that for every
$i\in\mathcal{M}$ and any $x\in\mathbb{R}^n$, we have
\[\|d\|\leq L \qquad \hbox{for all } d\in\partial f_i(x).\]
\label{ass.bounded_subgrad}
\end{assumption}

\subsection{Network Communication Model}

We define the {\it communication matrix} for the network at
time $k$ as $A(k)=[a_{ij}(k)]_{i,j\in\mathcal{M}}$. We assume a
probabilistic communication model, in which the sequence of
communication matrices $A(k)$ is assumed to be Markovian on the
{\it state variable}
$x(k)=[x_i(k)]_{i\in\mathcal{M}}\in\mathbb{R}^{n\times m}$.
Formally, let $\{n(k)\}_{k \in \mathbb{N}}$ be an independent
sequence of random variables defined in a probability space
$(\Omega, \mathcal{F},P) = \prod_{k=0}^\infty (\Omega',
\mathcal{F}',P')_k$, where $\{(\Omega', \mathcal{F}',P')_k\}_{k
\in \mathbb{N}}$ constitutes a sequence of identical
probability spaces. We assume there exists a function
$\psi:\mathbb{R}^{n\times m}\times\Omega'
\rightarrow\mathbb{R}^{m\times m}$ such that
\[
A(k)=\psi(x(k),n(k)).
\]
This Markovian communication model enables us to capture
settings where the agents' ability to communicate with each
other depends on their current estimates.

We assume there exists some underlying communication graph
$(\mathcal{M},\mathcal{E})$ that represents a `backbone' of the
network. That is, for each edge $e \in \mathcal{E}$, the two agents
linked by $e$ systematically attempt to communicate with each other
[see Eq.\ (\ref{eq.prob_model}) for the precise statement]. We do
not make assumptions on the communication (or lack thereof) between
agents that are not adjacent in $(\mathcal{M},\mathcal{E})$. We make
the following connectivity assumption on the graph
$(\mathcal{M},\mathcal{E})$.

\begin{assumption}[Connectivity]
The graph $(\mathcal{M},\mathcal{E})$ is strongly connected.
\label{connectivity}
\end{assumption}

The central feature of the model introduced in this paper is
that the probability of communication between two agents is
potentially small if their estimates are far apart. We
formalize this notion as follows: for all $(j,i) \in
\mathcal{E}$, all $k\geq 0$ and all $\overline{x}\in
\mathbb{R}^{m \times n}$,
\begin{equation}
P(a_{ij}(k)\geq\gamma |
x(k)=\overline{x}) \geq
\min\left\{\delta,\frac{K}{\|\overline{x}_i-\overline{x}_j\|^C}\right\}, \label{eq.prob_model}
\end{equation}
where $K$ and $C$ are real positive constants, and $\delta\in
(0,1]$. We included the parameter $\delta$ in the model to upper
bound the probability of communication when
$\|\overline{x}_i-\overline{x}_j\|^C$ is small. This model states
that, for any two nodes $i$ and $j$ with an edge between them, if
estimates $x_i(k)$ and $x_j(k)$ are close to each other, then there
is a probability at least $\delta$ that they communicate at time
$k$. However, if the two agents are far apart, the probability they
communicate can only be bounded by the inverse of a polynomial of
the distance between their estimates $\| x_i(k) - x_j(k)\|$. If the
estimates were to represent physical locations of wireless sensors,
then this bound would capture fading effects in the communication
channel.

We make two more technical assumptions to guarantee,
respectively, that the communication between the agents
preserves the average of the estimates, and the agents do not
discard their own information.

\begin{assumption}[Doubly Stochastic Weights]
The communication matrix $A(k)$ is doubly stochastic for all $k\geq
0$, i.e., for all $k\ge 0$, $a_{ij}(k)\geq 0$ for all
$i,j\in\mathcal{M}$, and $\sum_{i=1}^m a_{ij}(k)=1 $ for all
$j\in\mathcal{M}$ and $\sum_{j=1}^m a_{ij}(k)=1$ for all
$i\in\mathcal{M}$ with probability one. \label{ass.stoch_weights}
\end{assumption}

\begin{assumption}[Self Confidence]
There exists $\gamma>0$ such that $a_{ii}\geq\gamma$ for all agents
$i\in\mathcal{M}$ with probability one. \label{ass.self_conf}
\end{assumption}

The doubly stochasticity assumption on the matrices $A(k)$ is
satisfied when agents coordinate their weights when exchanging
information, so that $a_{ij}(k)=a_{ji}(k)$ for all
$i,j\in\mathcal{M}$ and $k\geq 0$.\footnote{This will be
achieved when agents exchange information about their estimates
and ``planned" weights simultaneously and set their actual
weights as the minimum of the planned weights; see
\cite{nedic-ozdaglar} where such a coordination scheme is
described in detail.} The self-confidence assumption states
that each agent gives a significant weight to its own estimate.

\subsection{A Counterexample}

In this subsection, we construct an example to demonstrate that the
algorithm defined in Eqs.\ (\ref{convex-comb})-(\ref{proj-error})
does not necessarily solve the optimization problem given in Eq.\
(\ref{optim-prob}). The following proposition shows that there exist
problem instances where Assumptions
\ref{ass.bounded_subgrad}-\ref{ass.self_conf} hold and $X_i = X$ for
all $i \in \mathcal{M}$, however the sequence of estimates $x_i(k)$
(and the sequence of function values $f(x_i(k))$) diverge for some
agent $i$ with probability one.

\begin{proposition}\label{prop:counter} Let Assumptions \ref{ass.bounded_subgrad}, \ref{connectivity}, \ref{ass.stoch_weights} and
\ref{ass.self_conf} hold and let $X_i = X$ for all $i \in
\mathcal{M}$. Let $\{x_i(k)\}$ be the sequences generated by
the algorithm (\ref{convex-comb})-(\ref{proj-error}). Let $C
> 1$ in Eq.\ (\ref{eq.prob_model}) and let
the stepsize be a constant value $\alpha$. Then, there does not
exist a bound $M(m,L,\alpha) < \infty$ such that
\[ \liminf_{k \to \infty}|f(x_i(k)) - f^*| \leq M(m,L,\alpha)\]  with probability 1, for all agents $i \in \mathcal{M}$.
\end{proposition}

\begin{proof}  Consider a network consisting of two agents solving a one-dimensional minimization problem.
 The first agent's objective function is $f_1(x) = -x$, while the second agent's objective function
 is $f_2(x) = 2x$. Both agents'
 feasible sets are equal to $X_1 = X_2 = [0,\infty)$. Let $x_1(0) \geq x_2(0) \geq 0$. The elements of the communication matrix are given by
$$ a_{1,2}(k) = a_{2,1}(k) = \left\{
   \begin{array}{ll}
     \gamma, & \hbox{with probability}\qquad \min\left\{\delta,\frac{1}{|x_1(k)-x_2(k)|^C}\right\}; \\
     0, & \hbox{with probability}\qquad 1- \min\left\{\delta,\frac{1}{|x_1(k)-x_2(k)|^C}\right\},
   \end{array}
 \right.$$ for some $\gamma \in (0,1/2]$ and $\delta \in [1/2,1)$.

The optimal solution set of this multi-agent
 optimization problem is the singleton $X^* = \{0\}$ and the optimal solution is $f^* =
 0$. We now prove that $\lim_{k \to \infty} x_1(k) = \infty$
 with probability 1 implying that $\lim_{k \to \infty}
 |f(x_1(k))-f^*| = \infty$.

From the iteration in Eq.\ (\ref{eq.update_rule}), we have that for
any $k$,
\begin{eqnarray}
x_1(k+1) &=& a_{1,1}(k)x_1(k) + a_{1,2}(k)x_2(k) + \alpha \label{eq:it counter1}\\
x_2(k+1) &=& \max\{0, a_{2,1}(k)x_1(k) + a_{2,2}(k)x_2(k) - 2\alpha\} \label{eq:it counter2}.
\end{eqnarray} We do not need to project $x_1(k+1)$ onto $X_1 = [0,\infty)$ because $x_1(k+1)$ is non-negative if $x_1(k)$ and $x_2(k)$ are both non-negative.
Note that since $\gamma \leq 1/2$, this iteration preserves
$x_1(k) \geq x_2(k) \geq 0$ for all $k \in \mathbb{N}$.

We now show that for any $k \in \mathbb{N}$ and any $x_1(k)
\geq x_2(k) \geq 0$, there is probability at least $\epsilon >
0$ that the two agents will never communicate again, i.e.,
\begin{equation} P(a_{1,2}(k') = a_{2,1}(k') = 0 \hbox{ for all }
k' \geq k| x(k)) \geq \epsilon > 0.\label{eq:inf
prod}\end{equation} If the agents do not communicate on periods
$k, k+1,...,k+j-1$ for some $j \geq 1$, then \begin{eqnarray*}
x_1(k+j) - x_2(k+j) &=& (x_1(k+j) - x_1(k)) + (x_1(k) - x_2(k))
+ (x_2(k) - x_2(k+j))\\ &\geq& \alpha j + 0 + 0,
\end{eqnarray*} from Eqs.\ (\ref{eq:it counter1}) and (\ref{eq:it counter2}) and the fact that
$x_1(k) \geq x_2(k)$. Therefore, the communication probability at
period $k+j$ can be bounded by
\[ P(a_{1,2}(k+j) =0| x(k), a_{1,2}(k') = 0 \hbox{ for all } k' \in \{k,...,k+j-1\}) \geq 1 - \min\{\delta,(\alpha j)^{-C}\}.\]
Applying this bound recursively for all $j \geq k$, we obtain
\begin{eqnarray*}&& P(a_{1,2}(k') = 0 \hbox{ for all } k' \geq k|
x(k))\\&&\qquad = \prod_{j=0}^\infty P(a_{1,2}(k+j) = 0 | x(k),
a_{1,2}(k') = 0 \hbox{ for all } k' \in \{k,...,k+j-1\})\\&& \qquad
\geq \prod_{j=0}^\infty \left(1 - \min\{\delta,(\alpha
j)^{-C}\}\right)\end{eqnarray*} for all $k$ and all $x_1(k) \geq
x_2(k)$. We now show that $\prod_{j=0}^\infty \left(1 -
\min\{\delta,(\alpha j)^{-C}\}\right) > 0$ if $C>1$. Define the
constant $\overline{K} = \left\lceil\frac{2^{1 \over
C}}{\alpha}\right\rceil$. Since $\delta \geq 1/2$, we have that
$(\alpha j)^{-C} \leq \delta$ for $j \geq \overline{K}$. Hence, we
can separate the infinite product into two components:
\[ \prod_{j=0}^\infty \left(1 - \min\{\delta,(\alpha
j)^{-C}\}\right) \geq \left[\prod_{j<\overline{K}} \left(1 -
\min\{\delta,(\alpha j)^{-C}\}\right)\right]\left[
\prod_{j\geq\overline{K}} \left(1 - (\alpha
j)^{-C}\right)\right].\] Note that the term in the first
brackets in the equation above is a product of a finite number
of strictly positive numbers and, therefore, is a strictly
positive number. We, thus, have to show only that
$\prod_{j\geq\overline{K}} \left(1 - (\alpha j)^{-C}\right)
> 0$. We can bound this product by
\begin{eqnarray*}&& \prod_{j\geq\overline{K}} \left(1 - (\alpha
j)^{-C}\right) = \exp\left( \log
\left(\prod_{j\geq\overline{K}} \left(1 - (\alpha
j)^{-C}\right)\right)\right)\\&&\qquad = \exp  \left(
\sum_{j\geq\overline{K}} \log \left(1 - (\alpha
j)^{-C}\right)\right)  \geq   \exp\left(
\sum_{j\geq\overline{K}} -(\alpha
j)^{-C}\log(4)\right),\end{eqnarray*} where the inequality
follows from $\log(x) \geq (x-1)\log(4)$ for all $x \in
[1/2,1]$. Since $C>1$, the sum $\sum_{j\geq\overline{K}}
(\alpha j)^{-C}$ is finite and $\prod_{j=0}^\infty \left(1 -
\min\{\delta,(\alpha j)^{-C}\}\right) > 0$, yielding Eq.\
(\ref{eq:inf prod}).

Let $K^*$ be the (random) set of periods when agents communicate,
i.e., $a_{1,2}(k) = a_{2,1}(k) = \gamma$ if and only if $k \in K^*$.
For any value $k \in K^*$ and any $x_1(k) \geq x_2(k)$, there is
probability at least $\epsilon$ that the agents do not communicate
after $k$. Conditionally on the state, this is an event independent
of the history of the algorithm by the Markov property. If $K^*$ has
infinitely many elements, then by the Borel-Cantelli Lemma we obtain
that, with probability 1, for infinitely many $k$'s in $K^*$ there
is no more communication between the agents after period $k$. This
contradicts the infinite cardinality of $K^*$. Hence, the two agents
only communicate finitely many times and $\lim_{k \to \infty} x_1(k)
= \infty$ with probability 1.
\end{proof}

The proposition above shows the algorithm given by Eqs.\
(\ref{convex-comb})-(\ref{proj-error}) does not, in general, solve
the global optimization problem (\ref{optim-prob}). However, there
are two important caveats when considering the implications of this
negative result. The first one is that the proposition only applies
if $C > 1$. We leave it is an open question whether the same
proposition would hold if $C \leq 1$. The second and more important
caveat is that we considered only a constant stepsize in Proposition
\ref{prop:counter}. The stepsize is typically a design choice and,
thus, could be chosen to be diminishing in $k$ rather than a
constant.
 In the subsequent sections, we prove that the
algorithm given by Eqs.\ (\ref{convex-comb})-(\ref{proj-error}) does
indeed solve the optimization problem of Eq.\ (\ref{optim-prob}),
under appropriate assumptions on the stepsize sequence.

\section{Analysis of Information
Exchange}\label{sec:stochastic}

\subsection{The Disagreement Metric}

In this section, we consider how some information that a given agent
$i$ obtains at time $s$ affects a different agent $j$'s estimate
$x_j(k)$ at a later time $k \geq s$. In particular, we introduce a
disagreement metric $\rho(k,s)$ that establishes how far some
information obtained by a given agent at time $s$ is from being
completely disseminated in the network at time $k$. The two
propositions at the end of this section provide bounds on
$\rho(k,s)$ under two different set of assumptions.

 In view of the the linear representation in
Eqs.\ (\ref{convex-comb})-(\ref{proj-error}), we can express
the evolution of the estimates using products of matrices: for
any $s\geq 0$ and any $k\geq s$, we define the {\it transition
matrices} as
\[
\Phi(k,s)= A(s)A(s+1)\cdots A(k-1)A(k)\qquad \hbox{for all $s$ and
$k$ with $k\ge s$}.
\]
Using the transition matrices, we can relate the estimates at
time $k$ to the estimates at time $s<k$ as follows: for all
$i$, and all $k$ and $s$ with $k>s$,
\begin{eqnarray}
x_i(k+1) = \sum_{j=1}^m[\Phi(k,s)]_{ij}x_j(s) &-&
\sum_{r=s+1}^k\sum_{j=1}^m [\Phi(k,r)]_{ij}\alpha(r-1)d_j(r-1) -
\alpha(k)d_i(k)\nonumber\\ & +& \sum_{r=s+1}^k \sum_{j=1}^m
[\Phi(k,r)]_{ij} e_j(r-1) + e_i(k).\label{evolution-est}
\end{eqnarray}

Observe from the iteration above that $[\Phi(k,s)]_{ij}$ determines
how the information agent $i$ obtains at period $s-1$ impacts agent
$j$'s estimate at period $k+1$. If $[\Phi(k,s)]_{ij} = 1/m$ for all
agents $j$, then the information agent $i$ obtained at period $s-1$
is evenly distributed in the network at time $k+1$. We, therefore,
introduce the {\it disagreement metric} $\rho$,
\begin{equation}\label{eq:def-over-rho} \rho(k,s) =
\max_{i,j\in\mathcal{M}}\left|[\Phi(k,s)]_{ij} -
\frac{1}{m}\right| \qquad \hbox{ for all } k\geq s\geq 0,
\end{equation}
which, when close to zero, establishes that all information obtained
at time $s-1$ by all agents is close to being evenly distributed in
the network by time $k+1$.

\subsection{Propagation of Information}

The analysis in the rest of this section is intended to produce
upper bounds on the disagreement metric $\rho(k,s)$. We start our
analysis by establishing an upper bound on the maximum distance
between estimates of any two agents at any time $k$. In view of our
communication model  [cf.\ Eq.\ (\ref{eq.prob_model})], this bound
will be essential in constructing positive probability events that
ensure information gets propagated across the agents in the network.

\begin{lemma}
Let Assumptions \ref{ass.bounded_subgrad} and
\ref{ass.stoch_weights} hold. Let $x_i(k)$ be generated by the
update rule in (\ref{eq.update_rule}). Then, we have the
following upper bound on the norm of the difference between the
agent estimates: for all $k\ge 0$,
\[
\max_{i,h\in\mathcal{M}}\|x_i(k)-x_h(k)\|\leq\Delta+2mL\sum_{r=0}^{k-1}\alpha(r)+
2\sum_{r=0}^{k-1} \sum_{j=1}^m \|e_j(r)\|,
\]where $\Delta=2m\max_{j\in\mathcal{M}}\|x_j(0)\|$, and $e_j(k)$
denotes the projection error. \label{lem.max_distance}
\end{lemma}

\begin{proof}
Letting $s=0$ in Eq.\ (\ref{evolution-est}) yields,
\begin{eqnarray*}
&&x_i(k) = \sum_{j=1}^m[\Phi(k-1,0)]_{ij}x_j(0)\\&&\qquad\qquad -
\sum_{r=1}^{k-1}\sum_{j=1}^m [\Phi(k-1,r)]_{ij}\alpha(r-1)d_j(r-1) -
\alpha(k-1)d_i(k-1)\nonumber\\&&\qquad\qquad  + \sum_{r=1}^{k-1} \sum_{j=1}^m
[\Phi(k-1,r)]_{ij} e_j(r-1) + e_i(k-1).\end{eqnarray*} Since the
matrices $A(k)$ are doubly stochastic with probability one for
all $k$ (cf.\ Assumption \ref{ass.stoch_weights}), it follows
that the transition matrices $\Phi(k,s)$ are doubly stochastic
for all $k\ge s\ge0$, implying that every entry
$[\Phi(k,s)]_{ij}$ belongs to $[0,1]$ with probability one.
Thus, for all $k$ we have,
\begin{eqnarray*}
\|x_i(k)\|\leq\sum_{j=1}^m\|x_j(0)\|&+&
\sum_{r=1}^{k-1}\sum_{j=1}^m\alpha(r-1)\|d_j(r-1)\| +
\alpha(k-1)\|d_i(k-1)\| \\&+& \sum_{r=1}^{k-1} \sum_{j=1}^m
\|e_j(r-1)\| + \|e_i(k-1)\|.
\end{eqnarray*}
Using the bound $L$ on the subgradients, this implies
\[\|x_i(k)\|\leq\sum_{j=1}^m\|x_j(0)\| +
\sum_{r=0}^{k-1}mL\alpha(r) + \sum_{r=0}^{k-1} \sum_{j=1}^m
\|e_j(r)\|.\] Finally, the fact that
$\|x_i(k)-x_h(k)\|\leq\|x_i(k)\|+\|x_h(k)\|$ for every
$i,h\in\mathcal{M}$, establishes the desired result.
\end{proof}

The lemma above establishes a bound on the distance between the
agents' estimates that depends on the projection errors $e_j$, which
are endogenously determined by the algorithm. However, if there
exists some $M
> 0$ such that $\|e_i(k)\| \leq M \alpha(k)$ for all $i\in
\mathcal{M}$ and all $k \geq 0$, then lemma above implies that, with
probability 1,
$\max_{i,h\in\mathcal{M}}\|x_i(k)-x_h(k)\|\leq\Delta+2m(L+M)\sum_{r=0}^{k-1}\alpha(r)$.
Under the assumption that such an $M$ exists, we define the
following set for each $k \in \mathbb{N}$,
\begin{equation}\label{eq:def-R-M}
R_M(k) = \left\{ x \in \mathbb{R}^{m \times  n}\ |\
\max_{i,h\in\mathcal{M}}\|x_i(k)-x_h(k)\|\leq\Delta+2m(L+M)\sum_{r=0}^{k-1}\alpha(r)\right\}.
\end{equation}
This set represents the set of agent states which can be reached
when the agents use the projected subgradient algorithm.

We next construct a sequence of events, denoted by $G(\cdot)$, whose
individual occurrence implies that information has been propagated
from one agent to all other agents, therefore, implying a
contraction of the disagreement metric $\rho$.

 We say a
link $(j,i)$ is {\it activated at time $k$} when
$a_{ij}(k)\geq\gamma$, and we denote by $\mathcal{E}(k)$ the
set of such edges, i.e.,
\[\mathcal{E}(k) = \{(j,i)\ |\ a_{ij}(k)\geq\gamma\}.\]
Here we construct an event in which the edges of the graphs
${\cal E}(k)$ are activated sequentially over time $k$, so that
information propagates from every agent to every other agent in
the network.

To define this event, we fix a node $w\in\mathcal{M}$ and
consider {\it two directed spanning trees} rooted at $w$ in the
graph $(\mathcal{M},{\cal E})$: an in-tree $T_{in,w}$ and an
out-tree $T_{out,w}$. In $T_{in,w}$ there exists a directed
path from every node $i\ne w$ to $w$; while in $T_{out,w}$,
there exists a directed path from $w$ to every node $i\ne w$.
The strongly connectivity assumption imposed on
$(\mathcal{M},{\cal E})$ guarantees that these spanning trees
exist and each contains $m-1$ edges (see \cite{LP}).

We order the edges of these spanning trees in a way such that
on any directed path from a node $i\ne w$ to node $w$, edges
are labeled in nondecreasing order. Let us represent the edges
of the two spanning trees with the order described above as
\begin{equation}T_{in,w}=\{e_1,e_2,\ldots,e_{n-1}\},\qquad
T_{out,w}=\{f_1,f_2,\ldots,f_{n-1}\}.\label{treeorder}\end{equation}
For the in-tree $T_{in,w}$, we pick an arbitrary leaf node and
label the adjacent edge as $e_1$; then we pick another leaf
node and label the adjacent edge as $e_2$; we repeat this until
all leaves are picked. We then delete the leaf nodes and the
adjacent edges from the spanning tree $T_{in,r}$, and repeat
the same process for the new tree. For the out-tree
$T_{out,w}$, we proceed as follows: pick a directed path from
node $w$ to an arbitrary leaf and sequentially label the edges
on that path from the root node $w$ to the leaf; we then
consider a directed path from node $w$ to another leaf and
label the unlabeled edges sequentially in the same fashion; we
continue until all directed paths to all the leaves are
exhausted.

For all $l=1,\ldots,m-1$, and any time $k\geq 0$, consider the
events
\begin{align}
& B_l(k) = \{\omega \in \Omega\ |\ a_{e_l}(k+l-1)\ge\gamma\},\label{eq.in_event}\\
& D_l(k) = \{\omega \in \Omega\ |\
a_{f_l}(k+(m-1)+l-1)\ge\gamma\}\label{eq.out_event},
\end{align}
 and define,
\begin{equation}
G(k) = \bigcap_{l=1}^{m-1} \Big(B_l(k)\cap D_l(k)\Big).
\label{eq.G_event}
\end{equation}
For all $l=1,\ldots,m-1$, $B_l(k)$ denotes the event that edge
$e_l\in T_{in,w}$ is activated at time $k+l-1$, while $D_l(k)$
denotes the event that edge $f_l\in T_{out,w}$ is activated at
time $k+(m-1)+l-1$. Hence, $G(k)$ denotes the event in which
each edge in the spanning trees $T_{in,w}$ and $T_{out,w}$ are
activated sequentially following time $k$, in the order given
in Eq.\ (\ref{treeorder}).

The following result establishes a bound on the probability of
occurrence of such a $G(\cdot)$ event. It states that the
probability of an event $G(\cdot)$ can be bounded as if the
link activations were independent and each link activation had
probability of occurring at least
\[\min\left\{\delta,\frac{K}{(\Delta+2m(L+M)\sum_{r=1}^{k+2m-3}\alpha(r))^C}\right\},\]
where the $k+2m-3$ follows from the fact that event $G(\cdot)$
is an intersection of $2(m-1)$ events occurring consecutively
starting at period $k$.

\begin{lemma}
Let Assumptions \ref{ass.bounded_subgrad}, \ref{connectivity} and
\ref{ass.stoch_weights} hold. Let $\Delta$ denote the constant
defined in Lemma \ref{lem.max_distance}. Moreover, assume that there
exists $M>0$ such that $\|e_i(k)\|\leq M\alpha(k)$ for all $i$ and
$k\geq 0$. Then,
\begin{enumerate}
\item[(a)] For all $s\in\mathbb{N}$, $k\geq s$, and any
    state $\overline x\in R_M(s)$,
\[
P(G(k)|x(s)=\overline x) \geq
\min\left\{\delta,\frac{K}{(\Delta+2m(L+M)\sum_{r=1}^{k+2m-3}\alpha(r))^C}\right\}^{2(m-1)}.
\]
\item[(b)] For all $k\geq 0$, $u\geq 1$, and any state
    $\overline x\in R_M(k)$,
\begin{align*}
P&\left(\bigcup_{l=0}^{u-1}
G(k+2(m-1)l)\Bigg|x(k)=\overline{x}\right)\\&\geq 1-\left(1-\min\left\{\delta,\frac{K}{(\Delta+2m(L+M)\sum_{r=1}^{k+2(m-1)u-1}\alpha(r))^C}\right\}^{2(m-1)}\right)^u.
\end{align*}
\end{enumerate}
\label{lem.bound_prob_G}
\end{lemma}
\begin{proof} (a) The proof is based on the fact that the communication matrices $A(k)$ are Markovian on the state $x(k)$, for all time $k\geq 0$. First, note that
\begin{align}
P(G(k)|x(s)=\overline{x}) &=
P\left(\bigcap_{l=1}^{m-1} \Big(B_l(k)\cap
D_l(k)\Big)\Bigg|x(s)=\overline{x}\right)\nonumber\\&=
P\left(\bigcap_{l=1}^{m-1}
B_l(k)\Bigg|x(s)=\overline{x}\right)P\left(\bigcap_{l=1}^{m-1}
D_l(k)\Bigg|\bigcap_{l=1}^{m-1}
B_l(k),x(s)=\overline{x}\right).\label{eq.prob_G_split}
\end{align} To simplify notation, let $W = 2m(L + M)$. We show that for all $k\geq s$,
\begin{equation}\label{eq:markov-G}
\inf_{\overline{x} \in R_M(s)} P\left(\bigcap_{l=1}^{m-1}
B_l(k)\Bigg|x(s)=\overline{x}\right)\geq
\min\left\{\delta,\frac{K}{(\Delta+ W\sum_{r=1}^{k+2m-3}\alpha(r))^C}\right\}^{(m-1)}.
\end{equation}
We skip the proof of the equivalent bound for the second term
in Eq.\ (\ref{eq.prob_G_split}) to avoid repetition. By
conditioning on $x(k)$ we obtain for all $k \geq s$,
\begin{eqnarray*}
&&\inf_{\overline{x} \in R_M(s)} P\left(\bigcap_{l=1}^{m-1}
B_l(k)\middle|x(s)=\overline{x}\right) =\\&& \qquad   \inf_{\overline{x} \in R_M(s)}\int_{x' \in \mathbb{R}^{m \times  n}} P\left(\bigcap_{l=1}^{m-1}
B_l(k)\middle|x(k)=x', x(s)=\overline{x}\right)dP(x(k)=x'|x(s)=\overline{x}).
\end{eqnarray*}
Using the Markov Property, we see that conditional on $x(s)$ can be
removed from the right-hand side probability above, since $x(k)$
already contains all relevant information with respect to
$\cap_{l=1}^{m-1} B_l(k)$. By the definition of $R_M(\cdot)$ [see
Eq.\ (\ref{eq:def-R-M})], if $x(s) \in R_M(s)$, then $x(k) \in
R_M(k)$ for all $k \geq s$ with probability 1. Therefore,
\begin{equation}\label{eq:tranfs-s-k}
\inf_{\overline{x} \in R_M(s)} P\left(\bigcap_{l=1}^{m-1}
B_l(k)\middle|x(s)=\overline{x}\right) \geq    \inf_{\overline{x} \in R_M(k)} P\left(\bigcap_{l=1}^{m-1}
B_l(k)\middle|x(k)=x'\right).
\end{equation}
 By the definition of $B_1(k)$,
\begin{eqnarray}
\label{eq.probC_event}&& \inf_{\overline{x} \in R_M(k)}P\left(\bigcap_{l=1}^{m-1}
B_l(k)\middle|x(k)=\overline{x}\right) =\\
&&\qquad \inf_{\overline{x} \in R_M(k)}P(a_{e_1}(k)\geq\gamma|x(s)=\overline{x})P\left(\bigcap_{l=2}^{m-1}
B_l(k)\middle|a_{e_1}(k)\geq\gamma,x(k)=\overline{x}\right).\nonumber
\end{eqnarray}
Define
\[
Q(k)=\min\left\{\delta,\frac{K}{\left(\Delta+W\sum_{r=1}^{k}\alpha(r)\right)^C}\right\},
\]
and note that, in view of the assumption imposed on the norm of
the projection errors and based on Lemma
\ref{lem.max_distance}, we get
\[
\max_{i,h\in\mathcal{M}}\|x_i(k)-x_h(k)\|\leq\Delta+W\sum_{r=0}^{k-1}\alpha(r).
\]
Hence, from Eq.\ (\ref{eq.prob_model}) we have
\begin{equation}
P(a_{ij}(k)\geq\gamma | x(k)=\overline{x})\geq Q(k).
\label{eq.bound_Pk}
\end{equation}
Thus, combining Eqs.\ (\ref{eq.probC_event}) and
(\ref{eq.bound_Pk}) we obtain,
\begin{equation}
\inf_{\overline{x} \in R_M(k)}P\left(\bigcap_{l=1}^{m-1}
B_l(k)\Bigg|x(k)=\overline{x}\right) \geq
Q(k)\inf_{\overline{x} \in R_M(k)}P\left(\bigcap_{l=2}^{m-1}
B_l(k)\Bigg|a_{e_1}(k)\geq\gamma,x(k)=\overline{x}\right).
\label{eq.Prob_split}
\end{equation}
By conditioning on the state $x(k+1)$, and repeating the use of
the Markov property and the definition of $R_M(k+1)$, we can
bound the right-hand side of the equation above,
\begin{eqnarray}
\nonumber&&\inf_{\overline{x} \in R_M(k)}P\left(\bigcap_{l=2}^{m-1}
B_l(k)\Bigg|a_{e_1}(k)\geq\gamma,x(k)=\overline{x}\right) \\\nonumber &=&\inf_{\overline{x} \in R_M(k)} \int_{x'} P\left(\bigcap_{l=2}^{m-1}
B_l(k)\middle|x(k+1)=x'\right)dP(x(k+1)=x'|a_{e_1}(k)\geq\gamma,x(k)=\overline{x})\\ &\geq&
\inf_{x' \in R_M(k+1)} P\left(\bigcap_{l=2}^{m-1}
B_l(k)\middle|x(k+1)=x'\right).\label{eq:markov-repeat}
\end{eqnarray}
Combining Eqs.\ (\ref{eq.probC_event}), (\ref{eq.Prob_split})
and (\ref{eq:markov-repeat}), we obtain
\[ \inf_{\overline{x} \in R_M(k)}P\left(\bigcap_{l=1}^{m-1}
B_l(k)\middle|x(k)=\overline{x}\right) \geq Q(k) \inf_{\overline{x}
\in R_M(k+1)} P\left(\bigcap_{l=2}^{m-1}
B_l(k)\middle|x(k+1)=x'\right).\] Repeating this process for all
$l=1,...,m-1$, this yields
\[ \inf_{\overline{x} \in R_M(k)}P\left(\bigcap_{l=1}^{m-1}
B_l(k)\middle|x(k)=\overline{x}\right) \geq \prod_{l=1}^{m-1} Q(k+l-1).\]
Since $Q$ is a decreasing function, $\prod_{l=1}^{m-1} Q(k+l-1)
\geq Q(k+2m-3)^{m-1}$. Combining with Eq.\
(\ref{eq:tranfs-s-k}), we have that for all $k \geq s$
\[ \inf_{\overline{x} \in R_M(s)}P\left(\bigcap_{l=1}^{m-1}
B_l(k)\middle|x(s)=\overline{x}\right) \geq Q(k+2m-3)^{m-1},\] producing the desired Eq.\ (\ref{eq:markov-G}).\\

(b) Let $G^c(k)$ represent the complement of $G(k)$. Note that
\begin{align*} P\left(\bigcup_{l=0}^{u-1}
G(k+2(m-1)l)\Bigg|x(k)=\overline{x}\right) = 1 -
P\left(\bigcap_{l=0}^{u-1}
G^c(k+2(m-1)l)\Bigg|x(k)=\overline{x}\right).
\end{align*} By conditioning on $G^c(k)$, we obtain
 \begin{eqnarray*} &&P\left(\bigcap_{l=0}^{u-1}
G^c(k+2(m-1)l)\Bigg|x(k)=\overline{x}\right) =\\&&\qquad P\left(G^c(k)\middle|x(k)=\overline{x}\right)P\left(\bigcap_{l=1}^{u-1}
G^c(k+2(m-1)l)\Bigg|G^c(k), x(k)=\overline{x}\right).\end{eqnarray*} We bound the term $P\left(G^c(k)\middle|x(k)=\overline{x}\right)$ using the result from part $(a)$. We bound the second term in the right-hand side of the equation above using the Markov property and the definition of $R_M(\cdot)$, which is the same technique from part $(a)$,
\begin{eqnarray*}
\nonumber&&\sup_{\overline{x} \in R_M(k)}P\left(\bigcap_{l=1}^{u-1}
G^c(k+2(m-1)l)\middle|G^c(k),x(k)=\overline{x}\right) \\\nonumber&&\qquad =\sup_{\overline{x} \in R_M(k)} \int_{x'} P\left(\bigcap_{l=1}^{u-1}
G^c(k+2(m-1)l)\middle|x(k+2(m-1))=x'\right)~\times\\&&\qquad\qquad\qquad\qquad\qquad dP(x(k+2(m-1))=x'|G^c(k),x(k)=\overline{x})\\\nonumber &&\qquad \leq
\sup_{\overline{x} \in R_M(k+2(m-1))} P\left(\bigcap_{l=1}^{u-1}
G^c(k+2(m-1)l)\middle|x(k+2(m-1))=x'\right).
\end{eqnarray*} The result follows by repeating the bound above $u$ times.
\end{proof}

The previous lemma bounded the probability of an event
$G(\cdot)$ occurring. The following lemma shows the implication
of the event $G(\cdot)$ for the disagreement metric.

\begin{lemma}
Let Assumptions \ref{connectivity}, \ref{ass.stoch_weights} and
\ref{ass.self_conf} hold. Let $t$ be a positive integer, and
let there be scalars $s<s_1<s_2<\cdots <s_t<k$, such that
$s_{i+1}-s_i\geq 2(m-1)$ for all $i=1,\dots,t-1$. For a fixed
realization $\omega \in \Omega$, suppose that events $G(s_i)$
occur for each $i=1,\dots,t$. Then,
\[
\rho(k,s)\leq
2\left(1+\frac{1}{\gamma^{2(m-1)}}\right)\left(1-\gamma^{2(m-1)}\right)^t.
\]
\label{lem.metric_bound}
\end{lemma} We skip the
proof of this lemma since it would mirror the proof of Lemma 6
in \cite{Ilanoptim}.

\subsection{Contraction Bounds}

In this subsection, we obtain two propositions that establish
 contraction bounds on the disagreement metric based on two
different sets of assumptions. For our first contraction bound,
we need the following assumption on the sequence of stepsizes.

\begin{assumption}\label{ass. small steps} \emph{(Limiting Stepsizes)}
The sequence of stepsizes $\{\alpha(k)\}_{k \in \mathbb{N}}$
satisfies
\[
\lim_{k\rightarrow\infty}k \log^p(k) \alpha(k)=0 \qquad \hbox{ for all $p < 1$}.
\]
\end{assumption}

The following lemma highlights two properties of stepsizes that
satisfy Assumption \ref{ass. small steps}: they are always
square summable and they are not necessarily summable. The
convergence results in Section \ref{sec:optim} require
stepsizes that are, at the same time, not summable and square
summable.

\begin{lemma}\label{stepsizeprop} Let $\{\alpha(k)\}_{k \in
\mathbb{N}}$ be a stepsize sequence that satisfies Assumption
\ref{ass. small steps}. Then, the stepsizes are square summable,
i.e, $\sum_{k=0}^\infty \alpha^2(k) < \infty$.
 Moreover, there exists a sequence of stepsizes $\{\overline{\alpha}(k)\}_{k \in \mathbb{N}}$ that
 satisfies Assumption \ref{ass. small steps} and is not
summable, i.e., $\sum_{k=0}^\infty \overline{\alpha}(k) = \infty$.
\end{lemma}

\begin{proof} From Assumption \ref{ass. small steps}, with
$p=0$, we obtain that there exists some $\overline{K} \in
\mathbb{N}$ such that $\alpha(k) \leq 1/k$ for all $k \geq
\overline{K}$. Therefore,
\[ \sum_{k=0}^\infty \alpha^2(k) \leq \sum_{k=0}^{\overline{K}-1} \alpha^2(k) +  \sum_{k=\overline{K}}^\infty
\frac{1}{k^2} \leq \overline{K}\max_{k \in
\{0,...,\overline{K}-1\}} \alpha^2(k) + \frac{\pi^2}{6} <
\infty.\] Hence, $\{\alpha(k)\}_{k \in \mathbb{N}}$ is square
summable. Now, let $\overline{\alpha}(k) =
\frac{1}{(k+2)\log(k+2)}$ for all $k \in \mathbb{N}$. This
sequence of stepsizes satisfies Assumption \ref{ass. small
steps} and is not summable since for all $K' \in \mathbb{N}$
\[ \sum_{k=0}^{K'} \overline{\alpha}(k) \geq \log(\log(K'+2))\]
and $\lim_{K' \to \infty} \log(\log(K'+2)) = \infty$.
\end{proof}

The following proposition is one of the central results in our
paper. It establishes, first, that for any fixed $s$, the
expected disagreement metric $E[\rho(k,s)|x(s)=\overline{x}]$
decays at a rate of $e^{\sqrt{k-s}}$ as $k$ goes to infinity.
Importantly, it also establishes that, as $s$ grows, the
contraction bound for a fixed distance $k-s$ decays slowly in
$s$. This slow decay is quantified by a function $\beta(s)$
that grows to infinity slower than the polynomial $s^q$ for any
$q > 0$.

\begin{proposition} Let Assumptions \ref{ass.bounded_subgrad}, \ref{connectivity}, \ref{ass.stoch_weights},
\ref{ass.self_conf}, and \ref{ass. small steps} hold. Assume also
that there exists some $M>0$ such that $\|e_i(k)\|\leq M\alpha(k)$
for all $i \in \mathcal{M}$ and $k \in \mathbb{N}$. Then, there
exists a scalar $\mu > 0$, an increasing function
$\beta(s):\mathbb{N}\to\mathbb{R}_+$ and a function
$S(q):\mathbb{N}\to\mathbb{N}$ such that
\begin{equation} \beta(s) \leq s^q \qquad \hbox{ for all $q > 0$ and all $s \geq S(q)$}\label{eq:def beta}
\end{equation}\begin{equation}
\hbox{and } \qquad E[\rho(k,s)|x(s)=\overline{x}] \leq \beta(s)
e^{-\mu\sqrt{k-s}}\quad\mbox{for all } k\geq s \geq 0, ~
\overline{x} \in R_M(s).
\label{eq:contraction}\end{equation}\label{prop:contraction}
\end{proposition}

\begin{proof} \emph{Part 1.} The first step of the proof is to define two functions, $g(k)$ and $w(k)$, that respectively
bound the sum of the stepsizes up to time $k$ and the inverse
of the probability of communication at time $k$, and prove some
limit properties of the functions $g(k)$ and $w(k)$ [see Eqs.\
(\ref{eq:limit-g}) and (\ref{eq:limit-w})]. Define $g(k):
\mathbb{R}_+ \to \mathbb{R}_+$ to be the linear interpolation
of $\sum_{r=0}^{\lfloor k\rfloor} \alpha(k)$, i.e,
\[g(k) = \sum_{r=0}^{\lfloor k \rfloor} \alpha(k) + (k - \lfloor k \rfloor) \alpha(k-\lfloor k \rfloor + 1).\]
Note that $g$ is differentiable everywhere except at integer
points and $g'(k) = \alpha(k-\lfloor k \rfloor + 1) =
\alpha(\lceil k \rceil)$ at $k \notin \mathbb{N}$. We thus
obtain from Assumption \ref{ass. small steps} that for all $p <
1$,
\begin{equation}\label{eq:limit-g}\lim_{k\rightarrow\infty, k \notin \mathbb{N}}k \log^p(k)
g'(k) = \lim_{k\rightarrow\infty}\lceil k \rceil \log^p(\lceil
k \rceil) \alpha(\lceil k \rceil) = 0.\end{equation} Define
$w(k)$ according to \begin{equation}\label{eq:def-w} w(k) =
\frac{(\Delta +
2m(L+M)g(k))^{2(m-1)C}}{K^{2(m-1)}},\end{equation} where
$\Delta = 2m \max_{j \in \mathcal{M}} \|x_j(0)\|$ and $K$ and
$C$ are parameters of the communication model [see Eq.\
(\ref{eq.prob_model})]. We now show that for any $p < 1$,
\begin{equation}\label{eq:limit-w}\lim_{k\rightarrow\infty, k \notin \mathbb{N}}k \log^p(k)
w'(k) = 0.\end{equation} If $\lim_{k \to \infty} w(k) <
\infty$, then the equation above holds immediately from Eq.\
(\ref{eq:limit-g}). Therefore, assume
 $\lim_{k \to \infty} w(k) = \infty$. By L'Hospital's Rule, for
 any $q > 0$,
\begin{equation}\label{eq:lhospital} \lim_{k \to \infty, k \notin \mathbb{N}}
\frac{w(k)}{\log^q(k)} = \frac{1}{q}\lim_{k\rightarrow\infty, k
\notin \mathbb{N}}\frac{k w'(k)}{\log^{q-1}(k)}.\end{equation} At
the same time, if we take $w(k)$ to the power $\frac{1}{2(m-1)C}$
before using L'Hospital's Rule, we obtain that for any $q > 0$,
\begin{eqnarray*} \lim_{k \to \infty, k \notin \mathbb{N}}
\left(\frac{w(k)}{\log^q(k)}\right)^{\frac{1}{2(m-1)C}}
&=&\frac{1}{K^{1/C}}\lim_{k\rightarrow\infty, k \notin
\mathbb{N}}\frac{\Delta + 2m(L+M)
g(k)}{\log^{\frac{q}{2(m-1)C}}(k)}\\&=&
\frac{4m(m-1)(L+M)C}{K^{1/C}q}\lim_{k\rightarrow\infty, k \notin
\mathbb{N}}\frac{k g'(k)}{\log^{\frac{q}{2(m-1)C}-1}(k)} =
0,\end{eqnarray*} where the last equality follows from Eq.\
(\ref{eq:limit-g}). From the equation above, we obtain
\begin{equation}\label{eq:limit-w2}\lim_{k \to \infty, k
\notin \mathbb{N}} \frac{w(k)}{\log^q(k)} = \left[\lim_{k \to
\infty, k \notin \mathbb{N}}
\left(\frac{w(k)}{\log^q(k)}\right)^{\frac{1}{2(m-1)C}}\right]^{2(m-1)C}
= 0,\end{equation} which combined with Eq.\ (\ref{eq:lhospital}), yields
the desired Eq.\ (\ref{eq:limit-w}) for any $p = 1-q < 1$.\\

\emph{Part 2.} The second step of the proof involves defining a
family of events $\{H_i(s)\}_{i,s}$ that occur with probability at
least $\phi
> 0$. We will later prove that an occurrence of $H_i(s)$ implies a
contraction of the distance between the estimates. Let $h_i(s) = i +
\lceil w(2s)\rceil$ for any $i, s \in \mathbb{N}$, where $w(\cdot)$
is defined in Eq.\ (\ref{eq:def-w}). We say the event $H_i(s)$
occurs if one out of a sequence of $G$-events [see definition in
Eq.\ (\ref{eq.G_event})] starting after $s$ occurs. In particular,
$H_i(s)$ is the union of $h_i(s)$ $G$-events and is defined as
follows,
\[ H_i(s) = \bigcup_{j=1}^{h_i(s)} G\left(s + 2(m-1)\left(j-1+\sum_{r=1}^{i-1} h_r(s)\right)\right) \qquad \hbox{ for all $i$, $s \in
\mathbb{N}$},\] where $\sum_{r=1}^0 (\cdot) = 0$.
See Figure 1 for a graphic representation of the $H_i(s)$
events.
\begin{figure}
\centering
\includegraphics[width=6.5in]{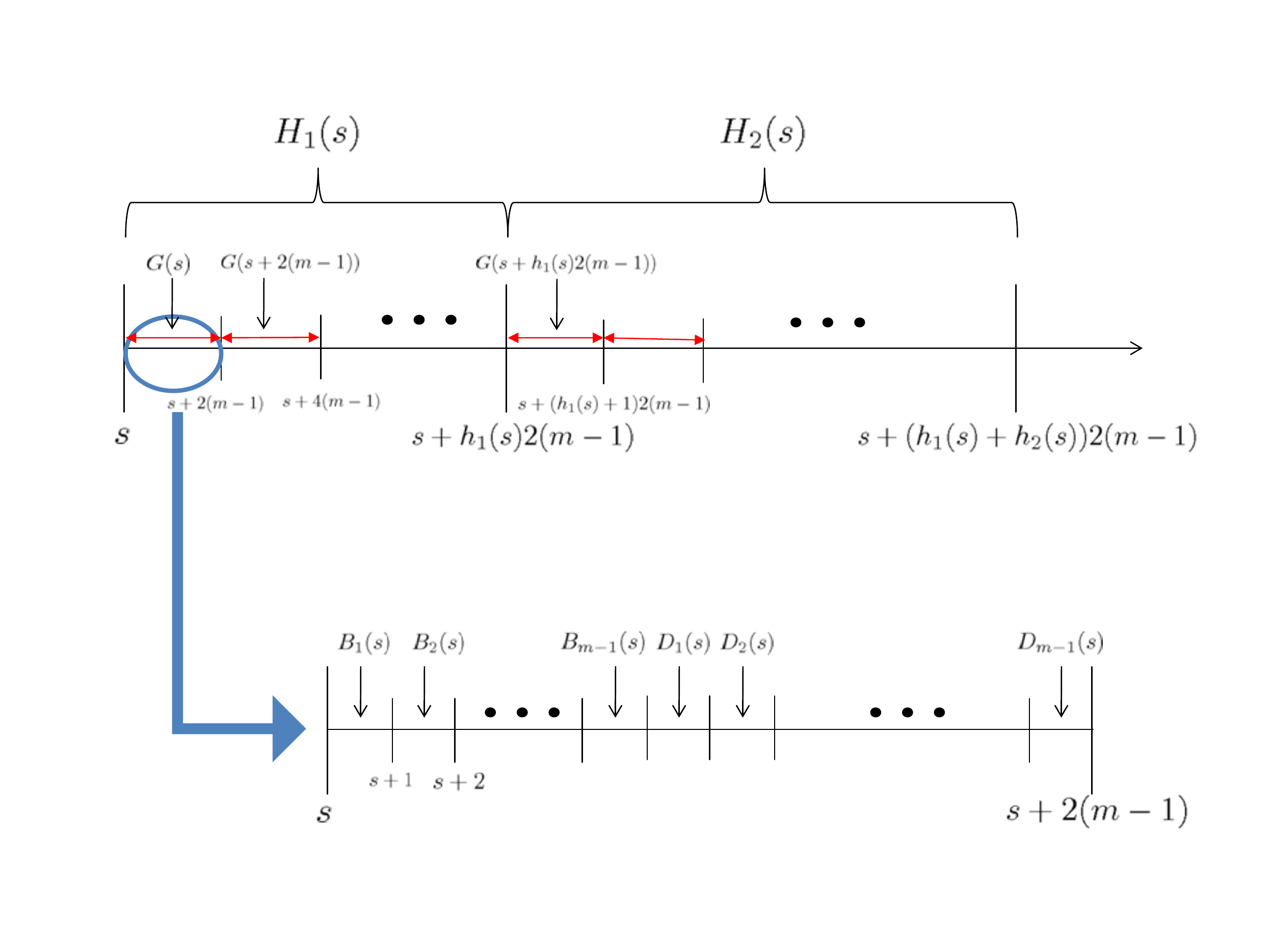}
\label{fig.events} \caption{The figure illustrates the three levels
of probabilistic events considered in the proof: the events $B_l(s)$
and $D_l(s)$, which represent the occurrence of communication over a
link (edge of the in-tree and out-tree, respectively); the events
$G(s)$ as defined in (\ref{eq.G_event}), with length $2(m-1)$ and
whose occurrence dictates the spread of information from any agent
to every other agent in the network; the events $H_i(s)$ constructed
as the union of an increasing number of events $G(s)$ so that their
probability of occurrence is guaranteed to be uniformly bounded away
from zero. The occurrence of an event $H_i(s)$ also implies the
spread of information from one agent to the entire network and, as a
result, leads to a contraction of the distance between the agents'
estimates.}
\end{figure}
We now show $P(H_i(s)|x(s) = \overline{x})$ for all $i,s \in
\mathbb{N}$ and all $\overline{x} \in R_M(s)$ [see definition
of $R_M(s)$ in Eq.\ (\ref{eq:def-R-M})]. From Lemma
\ref{lem.bound_prob_G}(a) and the definition of $w(\cdot)$, we
obtain that for all $\overline{x} \in R_M(s)$,
\[ P(G(s)|x(s)=\overline{x}) \geq
\min\left\{\delta^{2(m-1)},\frac{1}{w(s+2m-3)}\right\}.\]
Then, for all $s,i \in \mathbb{N}$ and all $\overline{x} \in
R_M(s)$,
\begin{eqnarray*}
P(H_i(s)|x(s)=\overline{x}) &=&
P\left(\bigcup_{j=1}^{h_i(s)} G\left(s + 2(m-1)\left(j-1+\sum_{r=1}^{i-1} h_r(s)\right)\right)\middle|x(s)=\overline{x}\right)\\
&\geq& 1 - \left(1 - \min\left\{\delta^{2(m-1)}, \frac{1}{w\left(s  + 2(m-1)\sum_{r=1}^{i} h_r(s)\right)}\right\}\right)^{h_i(s)},
\end{eqnarray*} where the inequality follows from Lemma \ref{lem.bound_prob_G}(b) and the fact that $w(\cdot)$ is a non-decreasing function.
Note that $h_r(s) \geq r$ for all $r$ and $s$, so that $s +
2(m-1)\sum_{r=1}^{i} h_r(s) \geq i^2$. Let $\hat{I}$ be the
smallest $i$ such that $w(i^2) \geq \delta^{-2(m-1)}$. We then
have that for all $i \geq \hat{I}$, all $s$ and all
$\overline{x} \in R_M(s)$,
\begin{eqnarray*}
P(H_i(s)|x(s)=\overline{x})
&\geq& 1 - \left(1 - \frac{1}{w\left(s  + 2(m-1)\sum_{r=1}^{i} h_r(s)\right)}\right)^{h_i(s)},
\end{eqnarray*}
Let $\tilde{I}$ be the maximum between $\hat{I}$ and the
smallest $i$ such that $w(i^2) > 1$. Using the inequality $(1 -
1/x)^x \leq e^{-1}$ for all $x \geq 1$, and multiplying and
dividing the exponent in the equation above by $w\left(s +
 2(m-1)\sum_{r=1}^{i} h_r(s)\right)$ we obtain
\begin{eqnarray*}
P(H_i(s)|x(s)=\overline{x}) \geq 1 - e^{-\frac{h_i(s)}{w\left(s +  2(m-1)\sum_{r=1}^{i} h_r(s)\right)}}
\end{eqnarray*} for all
$i \geq \tilde{I}$, all $s$ and all $\overline{x} \in R_M(s)$.
By bounding $h_r(s) \leq h_i(s)$ and replacing $h_i(s) = i +
\lceil w(2s)\rceil$, we obtain
\begin{eqnarray*}
P(H_i(s)|x(s)=\overline{x}) \geq 1 - e^{-\frac{i +
\lceil w(2s)\rceil}{w\left(s +  2(m-1)(i^2 + i
\lceil w(2s)\rceil)\right)}} \geq 1 - e^{-\frac{i +
 w(2s)}{w\left(s +  2(m-1)(i^2 + i w(2s) +i)\right)}}.
\end{eqnarray*} We now show there exists some $\overline{I}$ such that
\begin{equation}\label{eq:increasing} 1 - e^{-\frac{i +
 w(2s)}{w\left(s +  2(m-1)(i^2 + i
 w(2s) + i)\right)}} \qquad \hbox{is increasing in $i$ for all $i \geq \overline{I}$, $s \in \mathbb{N}$.}\end{equation}
The function above is increasing in $i$ if $\frac{i +
 w(2s)}{w\left(s +  2(m-1)(i^2 + i
 w(2s) + i)\right)}$ is increasing in $i$. The partial derivative of this function with respect to $i$ is positive if
\begin{eqnarray*} &&w\left(s +  2(m-1)(i^2 + i
 w(2s) + i)\right) - \\&& 2(m-1)(2i^2+i+3iw(2s)+w(2s)+w^2(2s))w'\left(s +  2(m-1)(i^2 + i
 w(2s) + i)\right) > 0\end{eqnarray*} at all points where the derivative $w'(\cdot)$ exists, that is, at non-integer values. If $i \geq \tilde{I}$, then $w\left(s +  2(m-1)(i^2 + i
 w(2s) + i)\right) > 1$ for all $s$ and it is thus sufficient to show
\[2(m-1)(2i^2+i+3iw(2s)+w(2s)+w^2(2s))w'\left(s +  2(m-1)(i^2 + i
 w(2s) + i)\right)  \leq 1\] in order to prove that Eq.\ (\ref{eq:increasing}) hold. The equation above holds if
\[2(m-1)(3i^2+4iw(2s)+w^2(2s))w'\left(2(m-1)(i^2 + i
 w(2s)) + s\right) \leq 1.\] From Eq.\ (\ref{eq:limit-w}) with $p=1/2$, we have that there exists some $N$ such that for all $x\geq N$, $w'(x) \leq \frac{1}{4x\sqrt{\log(x)}}$. For $i^2 \geq N$ and any $s \in \mathbb{N}$,
\begin{eqnarray*}&&2(m-1)(3i^2+4iw(2s)+w^2(2s))w'\left(2(m-1)(i^2 + i
 w(2s)) + s\right)\\&& \qquad \leq  \frac{2(m-1)(3i^2+4iw(2s)+w^2(2s))}{(2(m-1)(4i^2 +
 4 iw(2s)) + 4s)\sqrt{\log(2(m-1)(i^2 +
 iw(2s) + s)}}\\&&\qquad \leq \frac{3i^2+4iw(2s)+w^2(2s)}{4i^2 +
 4 iw(2s) + \frac{2}{m-1}s\sqrt{\log(i^2)}}.\end{eqnarray*} The term above is less than or equal to 1 if we select $i$ large enough such that $\frac{2}{m-1}s\sqrt{\log(i^2)} \geq w^2(2s)$ for all $s \in \mathbb{N}$ [see Eq.\ (\ref{eq:limit-w2}) with $q < 1/2$], thus proving there exists some $\overline{I}$ such that Eq.\ (\ref{eq:increasing}) holds. Hence, we obtain that for all $i, s \in \mathbb{N}$ and all $\overline{x} \in R_M(s)$,
\begin{eqnarray*}
&&P(H_i(s)|x(s)=\overline{x})\geq  \\&& \qquad \min_{j \in \{1,...,\overline{I}\}} \left\{ 1 - \left(1 - \min\left\{\delta^{2(m-1)}, \frac{1}{w\left(s  + 2(m-1)\sum_{r=1}^{j} h_r(s)\right)}\right\}\right)^{h_j(s)}\right\}.
\end{eqnarray*} Since $P(H_i(s)|x(s)=\overline{x}) > 0$ for all $i, s \in \mathbb{N}$ and all $\overline{x} \in R_M(s)$, to obtain the uniform lower bound on $P(H_i(s)|x(s)=\overline{x}) \geq \phi > 0$, it is sufficient to show that for all $i \in \{1,...,\overline{I}\}$ and all $\overline{x} \in R_M(s)$,
\[ \lim_{s \to \infty} P(H_i(s)|x(s)=\overline{x}) > 0.\] Repeating the steps above, but constraining $s$ to be large enough instead of $i$, we obtain there exists some $\tilde{S}$ such that for all $s \geq \tilde{S}$, all $i \in \mathbb{N}$ and $\overline{x} \in R_M(s)$,
\begin{eqnarray*}
P(H_i(s)|x(s)=\overline{x}) \geq 1 - e^{-\frac{i +
 w(2s)}{w\left(s +  2(m-1)(i^2 + i w(2s) +i)\right)}}.
\end{eqnarray*} Since there exists some $\hat{S}$ such that $w(2s) \leq \log(2s)$ for all $s \geq \hat{S}$ [see Eq.\ (\ref{eq:limit-w2}) with $q=1$], we obtain \begin{eqnarray*}
P(H_i(s)|x(s)=\overline{x}) \geq 1 - e^{-\frac{i +
 w(2s)}{w\left(s +  2(m-1)(i^2 + i \log(2s) +i)\right)}}.
\end{eqnarray*} for $s \geq \max\{\hat{S},\tilde{S}\}$ and all $i \in \mathbb{N}$ and $\overline{x} \in R_M(s)$. Note that for every $i$, there exists some $\overline{S}(i)$ such that for all $s  \geq \overline{S}(i)$  the numerator is greater than the denominator in the exponent above. Therefore, for all $i \in \mathbb{N}$ and $\overline{x} \in R_M(s)$,
\begin{eqnarray*}
\lim_{s \to \infty} P(H_i(s)|x(s)=\overline{x}) \geq 1 - e^{-1}.
\end{eqnarray*} Hence, there indeed exists some $\phi >0$ such that $P(H_i(s)|x(s)=\overline{x}) \geq \phi$ for all $i, s \in \mathbb{N}$ and $\overline{x} \in R_M(s)$.\\

\emph{Part 3.} In the previous step, we defined an event
$H_i(s)$ and proved it had probability at least $\phi > 0$ of
occurrence for any $i$ and $s$. We now determine a lower bound
on the number of possible $H$-events in an interval
$\{s,...,k\}$. The maximum number of possible $H$ events in the
interval $\{s,...,k\}$ is given by
\[ u(k,s) = \max\left\{t \in \mathbb{N}\ |\ s + 2(m-1)\sum_{i=1}^t h_i(s) \leq
k\right\}.\] Recall that $h_i(s) = i+\lceil w(2s)\rceil \leq i + w(2s) + 1$ to obtain
\[u(k,s) \geq \max\left\{t \in \mathbb{N}\ |\ \sum_{i=1}^t (i+ w(2s)+1) \leq
\frac{k-s}{2(m-1)}\right\}.\] By expanding the sum and adding $\left(\frac{3}{2}+w(2s)\right)^2$ to the left-hand side of the equation inside the maximization
 above, we obtain the following bound
\begin{eqnarray*} u(k,s) &\geq& \max\left\{t \in \mathbb{N}\ |\ t^2+ 3t +2w(2s)t +\left(\frac{3}{2}+w(2s)\right)^2 \leq \frac{k-s}{m-1}\right\}\\
&=& \max\left\{t \in \mathbb{N}\ |\ t+\frac{3}{2}+w(2s)\leq
\sqrt{\frac{k-s}{m-1}}\right\},
\end{eqnarray*} which yields the desired bound on $u(k,s)$,
\begin{equation}\label{eq:bound-u} u(k,s) \geq
\sqrt{\frac{k-s}{m-1}}-\frac{5}{2}-w(2s).\end{equation}

\emph{Part 4.} We now complete the proof of the proposition.
The following argument shows there is a high probability that
several $H$-events occur in a given $\{s,...,k\}$ interval and,
therefore, we obtain the desired contraction.

Let $I_i(s)$ be the indicator variable of the event $H_i(s)$,
that is $I_i(s) = 1$ if $H_i(s)$ occurs and $I_i(s) = 0$
otherwise.  For any $k \geq s \geq 0$, any $\overline{x} \in
R_M(s)$ and any $\delta >0$, the disagreement metric $\rho$
satisfies
\begin{eqnarray*}
&&E[\rho(k,s)|x(s)=\overline{x}] =\\&&\quad E\left[\rho(k,s)\middle|x(s)=\overline{x}, \sum_{i=1}^{u(k,s)} I_i(s) > \delta u(k,s)\right]
P\left(\sum_{i=1}^{u(k,s)} I_i(s) > \delta u(k,s)\middle|x(s)=\overline{x}\right)+\\&&\quad
E\left[\rho(k,s)\middle|x(s)=\overline{x}, \sum_{i=1}^{u(k,s)} I_i(s) \leq \delta u(k,s)\right]
P\left(\sum_{i=1}^{u(k,s)} I_i(s) \leq \delta u(k,s)\middle|x(s)=\overline{x}\right).
\end{eqnarray*} Since all the terms on the right-hand side of the equation above are less than or equal to 1, we obtain
\begin{eqnarray}
&&\label{eq:over-rho}E[\rho(k,s)|x(s)=\overline{x}] \leq\\&&\quad
E\left[\rho(k,s)\middle|x(s)=\overline{x}, \sum_{i=1}^{u(k,s)} I_i(s) > \delta u(k,s)\right]+
P\left(\sum_{i=1}^{u(k,s)} I_i(s) \leq \delta u(k,s)\middle|x(s)=\overline{x}\right).\nonumber
\end{eqnarray} We now bound the two terms in the right-hand side of Eq.\ (\ref{eq:over-rho}). Consider initially the first term. If $I_i(s) > \delta u(k,s)$,
then at least $\delta u(k,s)$ $H$-events occur, which by the
definition of $H_i(s)$ implies that at least $\delta u(k,s)$
$G$-events occur. From Lemma \ref{lem.metric_bound}, we obtain
\begin{equation}\label{eq:over-rho1}
E\left[\rho(k,s)\middle|x(s)=\overline{x}, \sum_{i=1}^{u(k,s)} I_i(s) > \delta u(k,s)\right]\leq
2\left(1+\frac{1}{\gamma^{2(m-1)}}\right)\left(1-\gamma^{2(m-1)}\right)^{\delta u(k,s)}
\end{equation} for all $\delta>0$.  We now consider the second term in the right-hand side of Eq.\ (\ref{eq:over-rho}).
The events $\{I_i(s)\}_{i=1,...,u(k,s)}$ all have probability at
least $\phi > 0$ conditional on any $x(s) \in R_M(s)$, but they are
not independent. However, given any $x(s+\sum_{i=1}^{j-1} h_j(r))
\in R(s+\sum_{i=1}^{j-1} h_j(r))$, the event $I_j(s)$ is independent
from the set of events $\{I_i(s)\}_{i=1,...,j-1}$ by the Markov
property. Therefore, we can define a sequence of independent
indicator variables $\{J_i(s)\}_{i=1,...,u(k,s)}$ such $P(J_i(s) =
1) = \phi$ and $J_i(s) \leq I_i(s)$ for all $i \in \{1,...,u(k,s)\}$
conditional on $x(s) \in R_M(s)$. Hence,
\begin{equation}\label{eq:over-rho2}P\left(\sum_{i=1}^{u(k,s)} I_i(s) \leq \delta
u(k,s)\middle|x(s)=\overline{x}\right) \leq
P\left(\sum_{i=1}^{u(k,s)} J_i(s) \leq \delta
u(k,s)\right),\end{equation} for any $\delta > 0$ and any
$\overline{x} \in R_M(s)$. By selecting $\delta =
\frac{\phi}{2}$ and using Hoeffding's Inequality, we obtain
\begin{equation}\label{eq:over-rho3}P\left(\frac{1}{u(k,s)}\sum_{i=1}^{u(k,s)} J_i(s) \leq
\frac{\phi}{2}\right) \leq
e^{-2\frac{\phi^2}{2^2}u(k,s)}.\end{equation} Plugging Eqs.\ (\ref{eq:over-rho1}), (\ref{eq:over-rho2})
 and (\ref{eq:over-rho3}), with $\delta = \phi/2$, into Eq.\
 (\ref{eq:over-rho}), we obtain
\begin{eqnarray*}
E[\rho(k,s)|x(s)=\overline{x}] \leq
2\left(1+\frac{1}{\gamma^{2(m-1)}}\right)\left(1-\gamma^{2(m-1)}\right)^{\frac{\phi}{2} u(k,s)} +
 e^{-\frac{\phi^2}{2}u(k,s)},
\end{eqnarray*} for all $k\geq s \geq 0$ and all $\overline{x} \in R_M(s)$. This implies there exists some  $\overline{\mu}_0, \overline{\mu}_1 > 0$ such
that $\rho(k,s)\leq \overline{\mu}_0 e^{-\overline{\mu}_1
u(k,s)}$ and, combined with Eq.\ (\ref{eq:bound-u}), we obtain
there exist some $K, \mu > 0$ such that
\[
E[\rho(k,s)|x(s)=\overline{x}] \leq
K e^{\mu\left(w(2s)-\sqrt{k-s}\right)} \qquad \hbox{ for all $k \geq s \geq 0, ~ x(s) \in R_M(s)$}.
\] Let $\beta(s) = Ke^{\mu w(2s)}$. Note that $\beta(\cdot)$ is an increasing function since $w(\cdot)$ is an increasing function.
To complete the proof we need to show that $\beta$ satisfies
condition stipulated in Eq.\ (\ref{eq:def beta}). From Eq.\
(\ref{eq:limit-w2}), with $q=1$, we obtain that
\[\lim_{s \to \infty, s
\notin \mathbb{N}} \frac{w(s)}{\log(s)} = 0.\] Note that since
$w(\cdot)$ is a continuous function, the limit above also
applies over the integers, i.e., $\lim_{s \to \infty}
\frac{w(s)}{\log(s)} = 0$.  Since $\lim_{s \to \infty}(s)$, for
any $q>0$, we have
\[0 = \lim_{s \to \infty} \frac{w(2s)}{\log(2s)} = \lim_{s \to \infty} \frac{\log(K) + \mu w(2s)}{q(-\log(2) + \log(2s))} =
\lim_{s \to \infty} \frac{\log(\beta(s))}{\log(s^q)}.
\] Let $S(q)$ be a scalar such that $\frac{\log(\beta(s))}{\log(s^q)} \leq 1$ for all $s \geq S(q)$. We thus obtain that $\beta(s) \leq s^q$ for all $s\geq S(q)$,
completing the proof of the proposition.
\end{proof}

The above proposition yields the desired contraction of the
disagreement metric $\rho$, but it assumes there exists some
$M>0$ such that $\|e_i(k)\|\leq M\alpha(k)$ for all $i \in
\mathcal{M}$ and $k \in \mathbb{N}$. In settings where we do
not have a guarantee that this assumption holds, we use the
proposition below. Proposition \ref{prop:compact-contraction}
instead requires that the sets $X_i$ be compact for each agent
$i$. With compact feasible sets, the contraction bound on the
disagreement metric follows not from the prior analysis in this
paper, but from the analysis of information exchange as if the
link activations were independent across time.

\begin{proposition} Let Assumptions \ref{connectivity}, \ref{ass.stoch_weights} and \ref{ass.self_conf} hold.
Assume also that the sets $X_i$ are compact for all $i \in
\mathcal{M}$. Then, there exist scalars $\kappa, \mu
> 0$ such that for all $\overline{x} \in \prod_{i \in
\mathcal{M}} X_i$,
\begin{equation}
E[\rho(k,s)|x(s)=\overline{x}] \leq \kappa
e^{-\mu(k-s)}\quad\mbox{for all } k\geq s \geq 0.
\label{eq:compact-contraction}\end{equation}\label{prop:compact-contraction}
\end{proposition}
\begin{proof} From Assumption \ref{connectivity}, we have that there exists a set of edges
 $\mathcal{E}$ of the strongly connected graph $(\mathcal{M}, \mathcal{E})$ such that for all
  $(j,i) \in \mathcal{E}$, all $k\geq 0$ and all $\overline{x}\in\mathbb{R}^{m \times n}$,
\[
P(a_{ij}(k)\geq\gamma | x(k)=\overline{x}) \geq
\min\left\{\delta,\frac{K}{\|\overline{x}_i-\overline{x}_j\|^C}\right\}.\]
The function
$\min\left\{\delta,\frac{K}{\|\overline{x}_i-\overline{x}_j\|^C}\right\}$
is continuous and, therefore, it attains its optimum when
minimized over the compact set $\prod_{i\in\mathcal{M}} X_i$,
i.e.,
\[\inf_{\overline{x} \in \prod_{i\in\mathcal{M}}
X_i}
\min\left\{\delta,\frac{K}{\|\overline{x}_i-\overline{x}_j\|^C}\right\}
= \min_{\overline{x} \in \prod_{i\in\mathcal{M}} X_i}
\min\left\{\delta,\frac{K}{\|\overline{x}_i-\overline{x}_j\|^C}\right\}.\]
Since the function
$\min\left\{\delta,\frac{K}{\|\overline{x}_i-\overline{x}_j\|^C}\right\}$
is strictly positive for any $\overline{x} \in \mathbb{R}^{m \times
n}$, we obtain that there exists some positive $\epsilon$ such that
\[\epsilon = \inf_{\overline{x} \in \prod_{i\in\mathcal{M}}
X_i}
\min\left\{\delta,\frac{K}{\|\overline{x}_i-\overline{x}_j\|^C}\right\}
> 0.\] Hence, for all $(j,i) \in \mathcal{E}$, all $k\geq 0$ and all
$\overline{x}\in\prod_{i\in\mathcal{M}} X_i$,
\begin{equation}\label{eq:rep-lemma7} P(a_{ij}(k)\geq\gamma |
x(k)=\overline{x}) \geq \epsilon.\end{equation} Since there is a
uniform bound on the probability of communication for any given edge
in $\mathcal{E}$ that is independent of the state $x(k)$, we can use
an extended version of Lemma 7 from \cite{Ilanoptim}. In particular,
Lemma 7 as stated in \cite{Ilanoptim} requires the communication
probability along edges to be independent of $x(k)$ which does not
apply here, however, it can be extended with straightforward
modifications to hold if the independence assumption were to be
replaced by the condition specified in Eq.\ (\ref{eq:rep-lemma7}),
implying the desired result.
\end{proof}

\section{Analysis of the Distributed Subgradient
Method}\label{sec:optim}

In this section, we study the convergence behavior of the agent
estimates $\{x_i(k)\}$ generated by the projected multi-agent
subgradient algorithm (\ref{eq.update_rule}). We first focus on
the case when the constraint sets of agents are the same, i.e.,
for all $i$, $X_i=X$ for some closed convex nonempty set. In
this case, we will prove almost sure consensus among agent
estimates and almost sure convergence of agent estimates to an
optimal solution when the stepsize sequence converges to 0
sufficiently fast (as stated in Assumption \ref{ass. small
steps}). We then consider the case when the constraint sets of
the agents $X_i$ are different convex compact sets and present
convergence results both in terms of almost sure consensus of
agent estimates and almost sure convergence of the agent
estimates to an optimal solution under weaker assumptions on
the stepsize sequence.

We first establish some key relations that hold under general
stepsize rules that are used in the analysis of both cases.

\subsection{Preliminary Relations}

The first relation measures the ``distance" of the agent
estimates to the intersection set $X=\cap_{i=1}^m X_i$. It will
be key in studying the convergence behavior of the projection
errors and the agent estimates. The properties of projection on
a closed convex set, subgradients, and doubly stochasticity of
agent weights play an important role in establishing this
relation.

\begin{lemma} Let Assumption \ref{ass.stoch_weights} hold.
Let $\{x_i(k)\}$ and $\{e_i(k)\}$ be the sequences generated by
the algorithm (\ref{convex-comb})-(\ref{proj-error}). For any
$z\in X=\cap_{i=1}^m X_i$, the following hold:
\begin{itemize}
\item[(a)] For all $k\ge 0$, we have
\begin{eqnarray*}\sum_{i=1}^m\|x_i(k+1)-z\|^2 &\le& \sum_{i=1}^m\|x_i(k)-z\|^2
+\a^2(k) \sum_{i=1}^m \|d_i(k)\|^2\\
&&\ -2\a(k)\sum_{i=1}^m\left(d_i(k)'(v_i(k)-z)\right)
 -\sum_{i=1}^m\|e_i(k)\|^2.\end{eqnarray*}
\item[(b)] Let also Assumption \ref{ass.bounded_subgrad}
    hold. For all $k\ge 0$, we have
\begin{equation}\sum_{i=1}^m\|x_i(k+1)-z\|^2 \le \sum_{i=1}^m\|x_i(k)-z\|^2
+\a^2(k)m L^2 -2\a(k) \sum_{i=1}^m
(f_i(v_i(k))-f_i(z)).\label{fvi-rel}\end{equation} Moreover, for all
$k\ge 0$, it also follows that
\begin{eqnarray}
\sum_{j=1}^m\|x_j(k+1)-z\|^2 &\le& \sum_{j=1}^m\|x_j(k)-z\|^2
+\a^2(k) mL^2 +2\a(k) L\sum_{j=1}^m \|x_j(k)- y(k)\|\cr
&&-2\a(k)\left(f(y(k))-f(z)\right),\label{fy-rel}
\end{eqnarray}
\end{itemize} \label{key-rel}
\end{lemma}

\begin{proof}
(a)\ Since  $x_i(k+1)=P_{X_i} [v_i(k)-\a(k) d_i(k)],$ it
follows from the property of the projection error $e_i(k)$ in
Eq.\ (\ref{projerror-feasdir}) that for any $z\in X,$
\[\|x_i(k+1)-z\|^2\le \|v_i(k)-\a(k) d_i(k)-z\|^2-\|e_i(k)\|^2.\]
By expanding the term $\|v_i(k)-\a(k) d_i(k)-z\|^2$, we obtain
\[\|v_i(k)-\a(k) d_i(k)-z\|^2=\|v_i(k)-z\|^2 +\a^2(k) \|d_i(k)\|^2
-2\a(k) d_i(k)'(v_i(k)-z).\] Since $v_i(k)=\sum_{j=1}^m a_{ij}(k)
x_j(k)$, using the convexity of the norm square function and
the stochasticity of the weights $a_{ij}(k)$, $j=1, \ldots,m$,
it follows that
\[\|v_i(k)-z\|^2\le \sum_{j=1}^m a_{ij}(k) \|x_j(k)-z\|^2.\]
Combining the preceding relations, we obtain
\begin{eqnarray*}
\|x_i(k+1)-z\|^2 &\le&  \sum_{j=1}^m a_{ij}(k)
\|x_j(k)-z\|^2+\a^2(k) \|d_i(k)\|^2\\
&& -2\a(k)d_i(k)'(v_i(k)-z) -\|e_i(k)\|^2.
\end{eqnarray*}
By summing the preceding relation over $i=1,\ldots,m,$ and
using the doubly stochasticity of the weights, i.e.,
\[\sum_{i=1}^m \sum_{j=1}^m a_{ij}(k) \|x_j(k)-z\|^2
=\sum_{j=1}^m\left(\sum_{i=1}^m a_{ij}(k)\right) \|x_j(k)-z\|^2
=\sum_{j=1}^m \|x_j(k)-z\|^2,\] we obtain the desired result.

\vskip .5pc

\noindent (b)\ Since $d_i(k)$ is a subgradient of $f_i(x)$ at
$x=v_i(k)$, we have
\[d_i(k)'(v_i(k)-z)\ge f_i(v_i(k))-f_i(z).\]
Combining this with the inequality in part (a), using
subgradient boundedness and dropping the nonpositive projection
error term on the right handside, we obtain
\[\sum_{i=1}^m\|x_i(k+1)-z\|^2 \le \sum_{i=1}^m\|x_i(k)-z\|^2
+\a^2(k)m L^2 -2\a(k) \sum_{i=1}^m (f_i(v_i(k))-f_i(z)),\] proving
the first claim. This relation implies that
\begin{eqnarray}
\sum_{j=1}^m\|x_j(k+1)-z\|^2 &\le& \sum_{j=1}^m\|x_j(k)-z\|^2
+\a^2(k) mL^2
-2\a(k)\sum_{i=1}^m\left(f_i(v_i(k))-f_i(y(k))\right)\cr
&&-2\a(k)\left(f(y(k))-f(z)\right).\qquad \label{eqn:mid}
\end{eqnarray}
In view of the subgradient boundedness and the stochasticity of
the weights, it follows
\[|f_i(v_i(k))-f_i(y(k))|\le L\|v_i(k)-y(k)\|
%=\left\|\sum_{j=1}^m a^i_j(k) x^j(k)-y(k)\right\|
\le L \sum_{j=1}^m a_{ij}(k) \|x_j(k)-y(k)\|,\] implying, by the
doubly stochasticity of the weights, that
\[\sum_{i=1}^m\left|f_i(v_i(k))-f_i(y(k))\right|\le
L \sum_{j=1}^m\left(\sum_{i=1}^m a_{ij}(k)\right) \|x_j(k)-y(k)\| =
L \sum_{j=1}^m \|x_j(k)-y(k)\|.\] By using this in relation
(\ref{eqn:mid}), we see that for any $z\in X$, and all $i$ and
$k,$
\begin{eqnarray*}
\sum_{j=1}^m\|x_j(k+1)-z\|^2 &\le& \sum_{j=1}^m\|x_j(k)-z\|^2
+\a^2(k) mL^2 +2\a(k) L\sum_{j=1}^m \|x_j(k)- y(k)\|\cr
&&-2\a(k)\left(f(y(k))-f(z)\right).
\end{eqnarray*}
\end{proof}

Our goal is to show that the agent disagreements
$\|x_i(k)-x_j(k)\|$ converge to zero. To measure the agent
disagreements $\|x_i(k)-x_j(k)\|$, we consider their average
$\frac{1}{m}\sum_{j=1}^m x_j(k)$, and consider the disagreement
of agent estimates with respect to this average. In particular,
we define
\begin{equation}y(k)=\frac{1}{m}\sum_{j=1}^m x_j(k)\qquad\hbox{for all }k.\label{yxaver}\end{equation} We
have
\[y(k+1)=\frac{1}{m}\sum_{i=1}^m v_i(k) -\frac{\a(k) }{m}\sum_{i=1}^m d_i(k)
+\frac{1}{m}\sum_{i=1}^m e_i(k).\] When the weights are doubly
stochastic, since $v_i(k)=\sum_{j=1}^m a_{ij}(k) x_j(k)$, it
follows that
\begin{equation}
y(k+1)=y(k) -\frac{\a(k) }{m}\sum_{i=1}^m d_i(k)
+\frac{1}{m}\sum_{i=1}^m e_i(k). \label{y_evol}
\end{equation}

Under our assumptions, the next lemma provides an upper bound
on the agent disagreements, measured by
$\Big\{\|x_i(k)-y(k)\|\Big\}$ for all $i$, in terms of the
subgradient bounds, projection errors and the disagreement
metric $\rho(k,s)$ defined in Eq.\ (\ref{eq:def-over-rho}).

\begin{lemma} Let Assumptions \ref{ass.bounded_subgrad} and \ref{ass.stoch_weights} hold.
Let $\{x_i(k)\}$ be the sequence generated by the algorithm
(\ref{convex-comb})-(\ref{proj-error}), and $\{y(k)\}$ be
defined in Eq.\ (\ref{y_evol}). Then, for all $i$ and $k\ge 2$,
an upper bound on $\|x_i(k)-y(k)\|$ is given by
\begin{eqnarray*}
\|x_i(k)-y(k)\| &\le & m\rho(k-1,0)\sum_{j=1}^m \|x_j(0)\| +
mL \sum_{r=0}^{k-2} \rho(k-1,r+1)\a(r) + 2\a(k-1) L\cr && +
\sum_{r=0}^{k-2} \rho(k-1,r+1)\sum_{j=1}^m \| e_j(r)\| +
\|e_i(k-1)\|+\frac{1}{m}\sum_{j=1}^m \|e_j(k-1)\|.
\end{eqnarray*}\label{est-xkyk}
\end{lemma}

\begin{proof} From Eq.\ (\ref{evolution-est}), we have for all $i$
and $k\ge s$,
\begin{eqnarray*}x_i(k+1) =
\sum_{j=1}^m[\Phi(k,s)]_{ij}x_j(s) &-& \sum_{r=s}^{k-1}\sum_{j=1}^m
[\Phi(k,r+1)]_{ij}\alpha(r)d_j(r) - \alpha(k)d_i(k)\nonumber\\ & +&
\sum_{r=s}^{k-1} \sum_{j=1}^m [\Phi(k,r+1)]_{ij} e_j(r) +
e_i(k).\end{eqnarray*} Similarly, using relation (\ref{y_evol}), we
can write for $y(k+1)$ and for all $k$ and $s$ with $k\ge s,$
\begin{eqnarray*}
y(k+1)= y(s)- \frac{1}{m}\sum_{r=s}^{k-1} \sum_{j=1}^m \a(r) d_j(r)
- \frac{\a(k)}{m}\sum_{i=1}^m  d_i(k)
 + \frac{1}{m}\sum_{r=s}^{k-1} \sum_{j=1}^m e_j(r)
+ \frac{1}{m}\sum_{j=1}^m e_j(k).
\end{eqnarray*}
Therefore, since $y(s)=\frac{1}{m} \sum_{j=1}^m x_j(s)$, we
have for $s=0,$
\begin{eqnarray*}
\|x_i(k)-y(k)\| &\le & \sum_{j=1}^m \left|[\Phi(k-1,0)]_{ij}
-\frac{1}{m}\right|\,\|x_j(0)\| \cr &&+ \sum_{r=0}^{k-2}
\sum_{j=1}^m \left|[\Phi(k-1,r+1)]_{ij}- \frac{1}{m}\right| \, \a(r)
\| d_j(r)\|\cr &&+ \a(k-1)\|d_i(k-1)\|
+\frac{\a(k-1)}{m}\sum_{j=1}^m \| d_j(k-1)\|\cr && +
\sum_{r=0}^{k-2} \sum_{j=1}^m
\left|[\Phi(k-1,r+1)]_{ij}-\frac{1}{m}\right| \| e_j(r)\| \cr &&+
\|e_i(k-1)\|+\frac{1}{m}\sum_{j=1}^m \|e_j(k-1)\|.
\end{eqnarray*}
Using the metric $\rho(k,s) = \max_{i,j \in {\cal M}}
\left|[\Phi(k,s)]_{ij}-{1\over m}\right|$ for $k\ge s\ge 0$
[cf.\ Eq.\ (\ref{eq:def-over-rho})], and the subgradient
boundedness, we obtain for all $i$ and $k\ge2,$
\begin{eqnarray*}
\|x_i(k)-y(k)\| &\le & m\rho(k-1,0)\sum_{j=1}^m \|x_j(0)\| +
mL \sum_{r=0}^{k-2} \rho(k-1,r+1)\a(r) + 2\a(k-1) L\cr && +
\sum_{r=0}^{k-2} \rho(k-1,r+1)\sum_{j=1}^m \| e_j(r)\| +
\|e_i(k-1)\|+\frac{1}{m}\sum_{j=1}^m \|e_j(k-1)\|,
\end{eqnarray*}
completing the proof.
\end{proof}

In proving our convergence results, we will often use the
following result on the infinite summability of products of
positive scalar sequences with certain properties. This result
was proven for geometric sequences in \cite{constconsoptim}.
Here we extend it for general summable sequences.

\begin{lemma} \label{lemma:seq}
Let $\{\beta_l\}$ and $\{\gamma_k\}$ be positive scalar
sequences, such that $\sum_{l=0}^{\infty}\beta_l<\infty$ and
$\lim_{k\to\infty}\gamma_k=0.$ Then,
\[\lim_{k\to\infty}\sum_{\ell=0}^k \beta_{k-\ell}\gamma_\ell=0.\]
In addition, if $\sum_{k=0}^\infty \gamma_k<\infty,$ then
\[\sum_{k=0}^\infty \sum_{\ell=0}^k \beta_{k-\ell}\gamma_\ell<\infty.\]
\end{lemma}

\begin{proof}
Let $\epsilon>0$ be arbitrary.  Since $\gamma_k\to0$, there is
an index $K$ such that $\gamma_k\le\epsilon$ for all $k\ge K$.
For all $k\ge K+1$, we have
\[\sum_{\ell=0}^k\beta_{k-\ell}\gamma_\ell=
\sum_{\ell=0}^K\beta_{k-\ell}\gamma_\ell +
\sum_{\ell=K+1}^k\beta_{k-\ell}\gamma_\ell \le \max_{0\le t\le
K}\gamma_t \sum_{\ell=0}^K\beta_{k-\ell}
+\epsilon\sum_{\ell=K+1}^k\beta_{k-\ell}.\]

Since $\sum_{l=0}^{\infty}\beta_l<\infty$, there exists $B>0$
such that
$\sum_{\ell=K+1}^k\beta_{k-\ell}=\sum_{\ell=0}^{k-K-1}\beta_{\ell}\le
B$ for all $k\ge K+1$. Moreover, since
$\sum_{\ell=0}^K\beta_{k-\ell}=\sum_{\ell=k-K}^{k}\beta_{\ell}$,
it follows that for all $k\ge K+1$,
\[\sum_{\ell=0}^k\beta_{k-\ell}\gamma_\ell
\le \max_{0\le t\le K}\gamma_t \sum_{\ell=k-K}^{k}\beta_{\ell}
+\epsilon B.\] Therefore, using $\sum_{l=0}^{\infty}\beta_l<\infty$,
we obtain
\[\limsup_{k\to\infty}\sum_{\ell=0}^k\beta_{k-\ell}\gamma_\ell
\le \epsilon B.\] Since $\epsilon$ is arbitrary, we conclude that
$\limsup_{k\to\infty}\sum_{\ell=0}^k\beta_{k-\ell}\gamma_\ell=0$,
implying
\[\lim_{k\to\infty}\sum_{\ell=0}^k\beta_{k-\ell}\gamma_\ell=0.\]

Suppose now $\sum_k\gamma_k<\infty$. Then, for any integer
$M\ge1$, we have
\[\sum_{k=0}^M \left(\sum_{\ell=0}^k\beta_{k-\ell}\gamma_\ell\right)
=\sum_{\ell=0}^M\gamma_\ell\sum_{t=0}^{M-\ell}\beta_t \le
\sum_{\ell=0}^M\gamma_\ell B,\] implying that
\[\sum_{k=0}^\infty \left(\sum_{\ell=0}^k\beta_{k-\ell}\gamma_\ell\right)
\le B\sum_{\ell=0}^\infty\gamma_\ell<\infty.\]
\end{proof}

\subsection{Convergence Analysis when $X_i=X$ for all
$i$}\label{sec:sameconst}

In this section, we study the case when agent constraint sets
$X_i$ are the same. We study the asymptotic behavior of the
agent estimates generated by the algorithm
(\ref{eq.update_rule}) using Assumption \ref{ass. small steps}
on the stepsize sequence.

The next assumption formalizes our condition on the constraint
sets.

\begin{assumption}
The constraint sets $X_i$ are the same, i.e., $X_i=X$ for a
closed convex set X.\label{sameconst}
\end{assumption}

We show first that under this assumption, we can provide an
upper bound on the norm of the projection error $\|e_i(k)\|$ as
a function of the stepsize $\a(k)$ for all $i$ and $k\ge 0$.

\begin{lemma} Let Assumptions \ref{ass.bounded_subgrad} and \ref{sameconst} hold.
Let $\{e_i(k)\}$ be the projection error defined by
(\ref{proj-error}). Then, for all $i$ and $k\ge 0$, the
$e_i(k)$ satisfy
\[\|e_i(k)\|\le 2 L \alpha(k).\]\label{bdprojerror}
\end{lemma}

\begin{proof}
Using the definition of projection error in Eq.\
(\ref{proj-error}), we have
\[e_i(k) = x_i(k+1) - v_i(k) + \alpha(k)d_i(k).\]
Taking the norms of both sides and using subgradient
boundedness, we obtain
\[\|e_i(k)\| \le \|x_i(k+1) - v_i(k)\| + \alpha(k)L.\]
Since $v_i(k) = \sum_{j=1}^m a_{ij}(k) x_j(k)$,  the weight
vector $a_i(k)$ is stochastic, and $x_j(k)\in X_j=X$ (cf.\
Assumption \ref{sameconst}), it follows that $v_i(k)\in X$ for
all $i$. Using the nonexpansive property of projection
operation [cf.\ Eq.\ (\ref{nonexpan})] in the preceding
relation, we obtain
\[\|e_i(k)\| \le \|v_i(k)-\alpha(k)d_i(k) - v_i(k)\| + \alpha(k)L \le 2\alpha(k)L,\]
completing the proof.
\end{proof}

This lemma shows that the projection errors are bounded by the
scaled stepsize sequence under Assumption \ref{sameconst}.
Using this fact and an additional assumption on the stepsize
sequence, we next show that the expected value of the sequences
$\{\|x_i(k)-y(k)\|\}$ converge to zero for all $i$, thus
establishing mean consensus among the agents in the limit. The
proof relies on the bound on the expected disagreement metric
$\rho(k,s)$ established in Proposition \ref{prop:contraction}.
The mean consensus result also immediately implies that the
agent estimates reach almost sure consensus along a particular
subsequence.

\begin{proposition} Let Assumptions \ref{ass.bounded_subgrad}, \ref{connectivity}, \ref{ass.stoch_weights}, \ref{ass.self_conf}, and \ref{sameconst} hold.
Assume also that the stepsize sequence $\{\a(k)\}$ satisfies
Assumption \ref{ass. small steps}. Let $\{x_i(k)\}$ be the
sequence generated by the algorithm
(\ref{convex-comb})-(\ref{proj-error}), and $\{y(k)\}$ be
defined in Eq.\ (\ref{y_evol}). Then, for all $i$, we have
\[\lim_{k\to \infty} E[\|x_i(k)-y(k)\|] = 0,\quad \hbox{and}\]\[\liminf_{k\to \infty}
\|x_i(k)-y(k)\| = 0\qquad \hbox{with probability one.}\]
\label{mean-consensus-sameset}
\end{proposition}

\begin{proof}
From Lemma \ref{est-xkyk}, we have the following for all $i$
and $k\ge 2$,
\begin{eqnarray*} \|x_i(k)-y(k)\| &\le &
m\rho(k-1,0)\sum_{j=1}^m \|x_j(0)\| + mL \sum_{r=0}^{k-2}
\rho(k-1,r+1)\a(r) + 2\a(k-1) L\cr && + \sum_{r=0}^{k-2}
\rho(k-1,r+1)\sum_{j=1}^m \| e_j(r)\| +
\|e_i(k-1)\|+\frac{1}{m}\sum_{j=1}^m \|e_j(k-1)\|.
\end{eqnarray*}
Using the upper bound on the projection error from Lemma
\ref{bdprojerror}, $\|e_i(k)\|\le 2\alpha(k)L$ for all $i$ and
$k$, this can be rewritten as
\begin{eqnarray}
\|x_i(k)-y(k)\| \le  m\rho(k-1,0)\sum_{j=1}^m \|x_j(0)\| &+&
3mL \sum_{r=0}^{k-2} \rho(k-1,r+1)\a(r)  \cr &+& 6\a(k-1)
L.\label{bdcont}
\end{eqnarray}
Under Assumption \ref{ass. small steps} on the stepsize
sequence, Proposition \ref{prop:contraction} implies the
following bound for the disagreement metric $\rho(k,s)$: for
all $k\ge s\ge 0$,
\[E[\rho(k,s)]\le \beta(s)e^{-\mu \sqrt{k-s}},\] where $\mu$ is a positive scalar and
$\beta(s)$ is an increasing sequence such that
\begin{equation}\beta(s)\leq s^q \qquad \hbox{ for all $q > 0$ and all $s
\geq S(q)$},\label{betabd}\end{equation} for some integer
$S(q)$, i.e., for all $q>0$, $\beta(s)$ is bounded by a
polynomial $s^q$ for sufficiently large $s$ (where the
threshold on $s$, $S(q)$, depends on $q$). Taking the
expectation in Eq.\ (\ref{bdcont}) and using the preceding
estimate on $\rho(k,s)$, we obtain
\begin{eqnarray*}
E[\|x_i(k)-y(k)\|] \le  m\beta(0)e^{-\mu\sqrt{k-1}}\sum_{j=1}^m
\|x_j(0)\| &+& 3mL \sum_{r=0}^{k-2} \beta(r+1)e^{-\mu
\sqrt{k-r-2}}\a(r)\\ &+& 6\a(k-1) L.
\end{eqnarray*}
We can bound $\beta(0)$ by $\beta(0)\le S(1)$ by using Eq.\
(\ref{betabd}) with $q=1$ and the fact that $\beta$ is an
increasing sequence. Therefore, by taking the limit superior in
the preceding relation and using  $\alpha(k)\to 0$ as $k\to
\infty$, we have for all $i$,
\begin{eqnarray*}
\limsup_{k\to \infty}E[\|x_i(k)-y(k)\|] \le 3mL \sum_{r=0}^{k-2}
\beta(r+1)e^{-\mu \sqrt{k-r-2}}\a(r).
\end{eqnarray*}
Finally, note that
\[\lim_{k\to \infty} \beta(k+1)\alpha(k)\le \lim_{k\to \infty} (k+1) \alpha(k) = 0,\]
where the  inequality holds by using Eq.\ (\ref{betabd}) with
$q=1$ and the equality holds by Assumption \ref{ass. small
steps} on the stepsize. Since we also have $\sum_{k=0}^\infty
e^{-\mu \sqrt{k}}<\infty$, Lemma \ref{lemma:seq} applies
implying that
\[\lim_{k\to \infty} \sum_{r=0}^{k-2}
\beta(r+1)e^{-\mu \sqrt{k-r-2}}\a(r)=0.\] Combining the preceding
relations, we have
\[\lim_{k\to \infty}E[\|x_i(k)-y(k)\|]=0.\]
Using Fatou's Lemma (which applies since the random variables
$\|y(k) - x_i(k)\|$ are nonnegative for all $i$ and $k$), we
obtain
\[ 0\le E\Big[\liminf_{k\to \infty}\|y(k) - x_i(k)\|\Big] \le \liminf_{k\to \infty} E[\|y(k) - x_i(k)\|] \le 0.\]
Thus, the nonnegative random variable $\liminf_{k\to
\infty}\|y(k) - x_i(k)\|$ has expectation 0, which implies that
\[\liminf_{k\to \infty}\|y(k) - x_i(k)\|=0\qquad \hbox{with probability one}.\]
\end{proof}

The preceding proposition shows that the agent estimates reach
a consensus in the expected sense. We next show that under
Assumption \ref{sameconst}, the agent estimates in fact
converge to an almost sure consensus in the limit. We rely on
the following standard convergence result for sequences of
random variables, which is an immediate consequence of the
supermartingale convergence theorem (see Bertsekas and
Tsitsiklis \cite{distbook}).

\begin{lemma}Consider a probability space $(\Omega, F,P)$ and
let $\{F(k)\}$ be an increasing sequence of $\sigma$-fields
contained in $F$. Let $\{V(k)\}$ and $\{Z(k)\}$ be sequences of
nonnegative random variables (with finite expectation) adapted
to $\{F(k)\}$ that satisfy
\[E[V(k+1)\ |\ F(k)]\le V(k) + Z(k),\]
\[\sum_{k=1}^\infty E[Z(k)] <\infty.\]
Then, $V(k)$ converges with probability one, as $k\to
\infty$.\label{stocapproxlemma}
\end{lemma}

\begin{proposition} Let Assumptions \ref{ass.bounded_subgrad}, \ref{connectivity}, \ref{ass.stoch_weights}, \ref{ass.self_conf}, and \ref{sameconst} hold.
Assume also that the stepsize sequence $\{\a(k)\}$ satisfies
Assumption \ref{ass. small steps}. Let $\{x_i(k)\}$ be the
sequence generated by the algorithm
(\ref{convex-comb})-(\ref{proj-error}), and $\{y(k)\}$ be
defined in Eq.\ (\ref{y_evol}). Then, for all $i$, we have:
\begin{itemize}
\item[(a)]\quad $\sum_{k=2}^\infty \a(k)\|x_i(k)-y(k)\|
    <\infty$ with probability one.
\item[(b)]\quad $\lim_{k\to \infty} \|x_i(k)-y(k)\| = 0$
    with probability one.
\end{itemize}
\label{as-consensus}
\end{proposition}

\begin{proof}
(a)\ Using the upper bound on the projection error from Lemma
\ref{bdprojerror}, $\|e_i(k)\|\le 2\alpha(k)L$ for all $i$ and
$k$, in Lemma \ref{est-xkyk}, we have for all $i$ and $k\ge 2$,
\begin{eqnarray*}
\|x_i(k)-y(k)\| \le  m\rho(k-1,0)\sum_{j=1}^m \|x_j(0)\| + 3mL
\sum_{r=0}^{k-2} \rho(k-1,r+1)\a(r) + 6\a(k-1) L.
\end{eqnarray*}
By multiplying this relation with $\alpha(k)$, we obtain
\begin{eqnarray*}
\alpha(k)\|x_i(k)-y(k)\| \le  m\a(k)\rho(k-1,0)\sum_{j=1}^m
\|x_j(0)\| &+& 3mL \sum_{r=0}^{k-2} \rho(k-1,r+1)\a(k)\a(r)\\
&+& 6\a(k)\a(k-1) L.
\end{eqnarray*}
Taking the expectation and using the estimate from Proposition
\ref{prop:contraction}, i.e.,
\[E[\rho(k,s)]\le \beta(s)e^{-\mu \sqrt{k-s}}\qquad \hbox{for all }k\ge s\ge 0,\] where
$\mu$ is a positive scalar and $\beta(s)$ is a increasing
sequence such that
\begin{equation}\beta(s)\leq s^q \qquad \hbox{ for all $q > 0$ and all $s
\geq S(q)$},\label{as-betabd}\end{equation} for some integer $S(q)$,
we have
\begin{eqnarray*}
E[\alpha(k)\|x_i(k)-y(k)\|] &\le&  m\a(k)\beta(0)e^{-\mu
\sqrt{k-1}}\sum_{j=1}^m \|x_j(0)\| \\&& + 3mL \sum_{r=0}^{k-2}
\beta(r+1)e^{-\mu \sqrt{k-r-2}}\a(k)\a(r)+ 6\a(k)\a(k-1) L.
\end{eqnarray*}
Let $\xi(r) = \beta(r+1) \alpha(r)$ for all $r\ge 0$. Using the
relations $\a(k)\xi(r)\le \a^2(k)+\xi^2(r)$ and
$2\a(k)\a(k-1)\le \a^2(k)+\a^2(k-1)$ for any $k$ and $r$, the
preceding implies that
\begin{eqnarray*}
E[\alpha(k)\|x_i(k)-y(k)\|] &\le&  m\a(k)\beta(0)e^{-\mu
\sqrt{k-1}}\sum_{j=1}^m
\|x_j(0)\| + 3mL\sum_{r=0}^{k-2} e^{-\mu \sqrt{k-r-2}}\xi^2(r)\\
&&\qquad + 3L\a^2(k) \Big(m \sum_{r=0}^{k-2}e^{-\mu \sqrt{k-r-2}}
+1\Big) +3\a^2(k-1) L.
\end{eqnarray*}
Summing over $k\ge 2$, we obtain
\begin{eqnarray*}
\sum_{k=2}^\infty E[\alpha(k)\|x_i(k)-y(k)\|] &\le& m\sum_{j=1}^m
\|x_j(0)\| \beta(0) \sum_{k=2}^\infty \a(k)e^{-\mu \sqrt{k-1}} \\ &&
\quad + 3L\sum_{k=2}^\infty \left(\Big(m \sum_{r=0}^{k-2} e^{-\mu
\sqrt{k-r-2}}+1\Big)\a^2(k) + \a^2(k-1)\right) \\ && \quad +
3mL\sum_{k=2}^\infty \sum_{r=0}^{k-2}e^{-\mu \sqrt{k-r-2}}\xi^2(r).
\end{eqnarray*}
We next show that the right handside of the above inequality is
finite: Since $\lim_{k\to \infty}\a(k)= 0$ (cf.\ Assumption
\ref{ass. small steps}), $\beta(0)$ is bounded, and $\sum_k
e^{-\mu \sqrt{k}}< \infty$, Lemma \ref{lemma:seq} implies that
the first term is bounded. The second term is bounded since
$\sum_k \a^2(k)<\infty$ by Assumption \ref{ass. small steps}
and Lemma \ref{stepsizeprop}. Since $\xi(r) =
\beta(r+1)\alpha(r)$, we have for some small $\epsilon>0$ and
all $r$ sufficiently large
\[\xi^2(r) = \beta^2(r+1)\alpha^2(r)\le (r+1)^{2/3} \alpha^2(r)\le (r+1)^{2/3}  \frac{\epsilon}{r^2},\]
where the first inequality follows using the estimate in Eq.\
(\ref{as-betabd}) with $q=1/3$ and the second inequality
follows from Assumption \ref{ass. small steps}. This implies
that $\sum_k \xi^2(k)<\infty$, which combined with Lemma
\ref{lemma:seq} implies that the third term is also bounded.
Hence, we have
\[\sum_{k=2}^\infty E[\a(k)\|x_i(k)-y(k)\|] <\infty.\]
By the monotone convergence theorem, this implies that
\[E\Big[\sum_{k=2}^\infty \alpha(k)\|y(k) - x_i(k)\|\Big]<\infty,\]
and therefore
\[\sum_{k=2}^\infty \alpha(k)\|y(k) - x_i(k)\|<\infty \qquad \hbox{with probability }1,\]
concluding the proof of this part.

\vskip .5pc

\noindent (b)\ Using the iterations (\ref{subgradient-step})
and (\ref{y_evol}), we obtain for all $k\ge 1$ and $i$,
\begin{eqnarray*}
y(k+1)-x_i(k+1) = \Big(y(k)-\sum_{j=1}^m a_{ij}(k) x_j(k)\Big)& -&
\alpha(k)\Big({1\over m} \sum_{j=1}^m d_j(k) -d_i(k)\Big)\\& +&
\Big({1\over m} \sum_{j=1}^m e_j(k) - e_i(k)\Big).
\end{eqnarray*}
By the stochasticity of the weights $a_{ij}(k)$ and the
subgradient boundedness, this implies that
\[\|y(k+1)-x_i(k+1)\| \le \sum_{j=1}^m a_{ij}(k)\|y(k)-x_j(k)\| + 2L\alpha(k) + {2\over m}\sum_{j=1}^m \|e_j(k)\|.\]
Using the bound on the projection error from Lemma
\ref{proj-error}, we can simplify this relation as
\[\|y(k+1)-x_i(k+1)\| \le \sum_{j=1}^m a_{ij}(k)\|y(k)-x_j(k)\| + 6L\alpha(k).\]
Taking the square of both sides and using the convexity of the
squared-norm function $\|\cdot\|^2$, this yields
\[\|y(k+1)-x_i(k+1)\|^2\le \sum_{j=1}^m a_{ij}(k)\|y(k)-x_j(k)\|^2 + 12L \alpha(k)
\sum_{j=1}^n a_{ij}(k)\|y(k)-x_j(k)\| + 36L^2 \alpha(k)^2.\] Summing
over all $i$ and using the doubly stochasticity of the weights
$a_{ij}(k)$, we have for all $k\ge 1$,
\[\sum_{i=1}^m\|y(k+1)-x_i(k+1)\|^2\le \sum_{i=1}^m \|y(k)-x_i(k)\|^2 + 12L \alpha(k)
\sum_{i=1}^m \|y(k)-x_i(k)\| + 36L^2 m \alpha(k)^2.\] By part (a) of
this lemma, we have $\sum_{k=1}^\infty \alpha(k)\|y(k) -
x_i(k)\|<\infty$ with probability one. Since, we also have
$\sum_k \alpha^2(k) <\infty$ (cf.\ Lemma \ref{stepsizeprop}),
Lemma \ref{stocapproxlemma} applies and implies that
$\sum_{i=1}^m\|y(k)-x_i(k)\|^2$ converges with probability one,
as $k\to \infty$.

By Proposition \ref{mean-consensus-sameset}, we have
\[\liminf_{k\to \infty} \|x_i(k)-y(k)\| = 0\qquad \hbox{with probability one.}\]
Since $\sum_{i=1}^m\|y(k)-x_i(k)\|^2$ converges with
probability one, this implies that for all $i$,
\[\lim_{k\to \infty}\|x_i(k)-y(k)\|=0\qquad \hbox{with probability one},\]
completing the proof.
\end{proof}

We next present our main convergence result under Assumption
\ref{ass. small steps} on the stepsize and Assumption
\ref{sameconst} on the constraint sets.

\begin{theorem} Let Assumptions \ref{ass.bounded_subgrad}, \ref{connectivity}, \ref{ass.stoch_weights},
\ref{ass.self_conf} and \ref{sameconst} hold. Assume also that
the stepsize sequence $\{\a(k)\}$ satisfies $\sum_{k=0}^\infty
\alpha(k)=\infty$ and Assumption \ref{ass. small steps}. Let
$\{x_i(k)\}$ be the sequence generated by the algorithm
(\ref{convex-comb})-(\ref{proj-error}). Then, there exists an
optimal solution $x^*\in X^*$ such that for all $i$
\[\lim_{k\to \infty} x_i(k)  = x^*\qquad \hbox{with probability one}.\]
\end{theorem}

\begin{proof}
From Lemma \ref{key-rel}(b), we have for some $z^*\in X^*$
(i.e., $f(z^*)=f^*$),
\begin{eqnarray}
\sum_{j=1}^m\|x_j(k+1)-z^*\|^2 &\le& \sum_{j=1}^m\|x_j(k)-z^*\|^2
+\a^2(k) mL^2 +2\a(k) L\sum_{j=1}^m \|x_j(k)- y(k)\|\cr
&&-2\a(k)\left(f(y(k))-f^*\right),\label{iterates-rel}
\end{eqnarray}
[see Eq.\ (\ref{fy-rel})]. Rearranging the terms and summing
these relations over $k=0,\ldots,K$, we obtain
\begin{eqnarray*}
2 \sum_{k=0}^K \a(k)\left(f(y(k))-f^*\right) &\le&
\sum_{j=1}^m\|x_j(0)-z^*\|^2 - \sum_{j=1}^m\|x_j(K+1)-z^*\|^2 \\&&\
+mL^2 \sum_{k=0}^K \a^2(k)  +2L \sum_{k=0}^K \a(k)\sum_{j=1}^m
\|x_j(k)- y(k)\|.
\end{eqnarray*}
By letting $K\to \infty$ in this relation and using
$\sum_{k=0}^\infty \alpha^2(k)<\infty$ (cf.\ Lemma
\ref{stepsizeprop}) and $\sum_{k=0}^\infty \a(k) \sum_{j=1}^m
\|x_j(k)-y(k)\|<\infty$ with probability one, we obtain
\[ \sum_{k=0}^K \a(k)\left(f(y(k))-f^*\right) <\infty\qquad \hbox{with probability one}.\]
Since $x_i(k)\in X$ for all $i$, we have $y(k)\in X$ [cf.\ Eq.\
(\ref{yxaver})] and therefore $f(y(k))\ge f^*$ for all $k$. Combined
with the assumption $\sum_{k=0}^\infty \a(k)=\infty$, the preceding
relation implies
\begin{equation}\liminf_{k\to \infty}
f(y(k))=f^*.\label{yvalue}\end{equation}

By dropping the nonnegative term
$2\a(k)\left(f(y(k))-f^*\right)$ in Eq.\ (\ref{iterates-rel}),
we have
\begin{eqnarray}
\sum_{j=1}^m\|x_j(k+1)-z^*\|^2 \le\sum_{j=1}^m\|x_j(k)-z^*\|^2
+\a^2(k) mL^2 +2\a(k) L\sum_{j=1}^m \|x_j(k)- y(k)\|.
\end{eqnarray}
Since $\sum_{k=0}^\infty \a^2(k)<\infty$ and $\sum_{k=0}^\infty
\a(k) \sum_{j=1}^m \|x_j(k)-y(k)\|<\infty$ with probability
one, Lemma \ref{stocapproxlemma} applies and implies that
$\sum_{j=1}^m\|x_j(k)-z^*\|^2$ is a convergent sequence with
probability one for all $z^*\in X^*$. By Lemma
\ref{as-consensus}(b), we have $\lim_{k\to \infty}
\|x_i(k)-y(k)\|=0$ with probability one, therefore it also
follows that the sequence $\|y(k)-z^*\|$ is also convergent.
Since $y(k)$ is bounded, it must have a limit point. By Eq.\
(\ref{yvalue}) and the continuity of $f$ (due to convexity of
$f$ over $\mathbb{R}^n$), this implies that one of the limit
points
 of $\{y(k)\}$ must belong to $X^*$; denote this limit
point by $x^*$. Since the sequence $\{\|y(k)-x^*\|\}$ is
convergent, it follows that $y(k)$ can have a unique limit
point, i.e., $\lim_{k\to \infty} y(k)=x^*$ with probability
one. This and $\lim_{k\to\infty}\|x_i(k)-y(k)\|=0$ with
probability one imply that each of the sequences $\{x_i(k)\}$
converges to the same $x^*\in X^*$ with probability one.
\end{proof}

%%%%%%%%%%%%%%%%%%
%%%%%%%%%%%%%%%%%%%

\subsection{Convergence Analysis for Different Constraint
Sets}\label{sec:difconst}

In this section, we provide our convergence analysis for the
case when all the constraint sets $X_i$ are different. We show
that even when the constraint sets of the agents are different,
the agent estimates converge almost surely to an optimal
solution of problem (\ref{optim-prob}) under some conditions.
In particular, we adopt the following assumption on the
constraint sets.

\begin{assumption} \label{compactconst}
For each $i$, the constraint set $X_i$ is a convex and compact
set.
\end{assumption}

An important implication of the preceding assumption is that
for each $i$, the subgradients of the function $f_i$ at all
points $x\in X_i$ are uniformly bounded, i.e., there exists
some scalar $L>0$ such that for all $i$,
\[\|d\| \le L \qquad \hbox{for all }d\in \partial f_i(x)\hbox{ and all }x\in X_i.\]

Our first lemma shows that with different constraint sets and a
stepsize that goes to zero, the projection error $e_i(k)$
converges to zero for all $i$ along all sample paths.

\begin{lemma} Let Assumptions  \ref{ass.stoch_weights} and \ref{compactconst}
hold. Let $\{x_i(k)\}$ and $\{e_i(k)\}$ be the sequences
generated by the algorithm
(\ref{convex-comb})-(\ref{proj-error}). Assume that the
stepsize sequence satisfies $\alpha(k)\to 0$ as $k$ goes to
infinity.
\begin{itemize}
\item[(a)] For any $z\in X$, the scalar sequence
    $\sum_{i=1}^m \|x_i(k)-z\|^2$ is convergent.
\item[(b)] The projection errors $e_i(k)$ converge to zero
    as $k\to \infty$, i.e.,
\[\lim_{k\to \infty} \|e_i(k)\|=0\qquad \hbox{for all }i.\]
\end{itemize}
\label{error-behav}
\end{lemma}

\begin{proof}
(a)\ Using subgradient boundedness and the relation
$|d_i(k)'(v_i(k)-z)|\le \|d_i(k)\|\|v_i(k)-z\|$ in part (a) of
Lemma \ref{key-rel}, we obtain
\[\sum_{i=1}^m\|x_i(k+1)-z\|^2 \le \sum_{i=1}^m\|x_i(k)-z\|^2
+\a^2(k) mL^2 +2\a(k)L \sum_{i=1}^m \|v_i(k)-z\|
 -\sum_{i=1}^m\|e_i(k)\|^2.\]
Since $v_i(k)=\sum_{j=1}^m a_{ij}(k) x_j(k)$, using doubly
stochasticity of the weights, we have $\sum_{i=1}^m
\|v_i(k)-z\|\le \sum_{i=1}^m \|x_i(k)-z\|$, which when combined
with the preceding yields for any $z\in X$ and all $k\ge 0$,
\begin{equation}\sum_{i=1}^m\|x_i(k+1)-z\|^2 \le \sum_{i=1}^m\|x_i(k)-z\|^2
+\a^2(k) mL^2 +2\a(k)L \sum_{i=1}^m \|x_i(k)-z\|
-\sum_{i=1}^m\|e_i(k)\|^2.\label{reduced-rel}\end{equation}

Since $x_i(k)\in X_i$ for all $i$ and $X_i$ is compact (cf.\
Assumption \ref{compactconst}), it follows that the sequence
$\{x_i(k)\}$ is bounded for all $i$, and therefore the sequence
$\sum_{i=1}^m \|x_i(k)-z\|$ is bounded. Since $\a(k)\to 0$ as
$k\to \infty$, by dropping the nonnegative term
$\sum_{i=1}^m\|e_i(k)\|^2$ in Eq.\ (\ref{reduced-rel}), it
follows that
\begin{eqnarray*}
\limsup_{k\to \infty}\sum_{i=1}^m\|x_i(k+1)-z\|^2 &\le&
\liminf_{k\to \infty}\sum_{i=1}^m\|x_i(k)-z\|^2\\ &&\ + \lim_{k\to
\infty}\left(\a^2(k) mL^2 +2\a(k)L \sum_{i=1}^m \|x_i(k)-z\|\right)\\
& =& \liminf_{k\to \infty}\sum_{i=1}^m\|x_i(k)-z\|^2.
\end{eqnarray*}
Since the sequence $\sum_{i=1}^m \|x_i(k)-z\|^2$ is bounded,
the preceding relation implies that the scalar sequence
$\sum_{i=1}^m \|x_i(k)-z\|^2$ is convergent.

\vskip .5pc

\noindent (b)\ From Eq.\ (\ref{reduced-rel}), for any $z\in X$,
we have
\[\sum_{i=1}^m \|e_i(k)\|^2 \le
\sum_{i=1}^m\|x_i(k)-z\|^2 - \sum_{i=1}^m\|x_i(k+1)-z\|^2  +\a^2(k)
mL^2 +2\a(k)L \sum_{i=1}^m \|x_i(k)-z\|.\] Taking the limit superior
as $k\to \infty$, we obtain
\begin{eqnarray*}
\limsup_{k\to \infty}\sum_{i=1}^m \|e_i(k)\|^2 &\le& \lim_{k\to
\infty }\left(\sum_{i=1}^m\|x_i(k)-z\|^2 -
\sum_{i=1}^m\|x_i(k+1)-z\|^2\right)\\ && +\lim_{k\to
\infty}\Big(\a^2(k) mL^2 +2\a(k)L \sum_{i=1}^m
\|x_i(k)-z\|\Big),\end{eqnarray*} where the first term on the right
handside is equal to zero by the convergence of the sequence
$\sum_{i=1}^m\|x_i(k)-z\|^2$, and the second term is equal to
zero by $\lim_{k\to \infty}\a(k)=0$ and the boundedness of the
sequence $\sum_{i=1}^m \|x_i(k)-z\|$, completing the proof.
\end{proof}

The preceding lemma shows the interesting result that the
projection errors $\|e_i(k)\|$ converge to zero along all
sample paths even when the agents have different constraint
sets under the compactness conditions of Assumption
\ref{compactconst}. Similar to the case with $X_i=X$ for all
$i$, we next establish mean consensus among the agent
estimates. The proof relies on the convergence of projection
errors to zero and the bound on the disagreement metric
$\rho(k,s)$ from Proposition \ref{prop:compact-contraction}.
Note that this result holds for all stepsizes $\a(k)$ with
$\a(k)\to 0$ as $k\to \infty$.

\begin{proposition} Let Assumptions \ref{connectivity}, \ref{ass.stoch_weights},
\ref{ass.self_conf} and \ref{compactconst} hold. Let
$\{x_i(k)\}$ be the sequence generated by the algorithm
(\ref{convex-comb})-(\ref{proj-error}), and $\{y(k)\}$ be
defined in Eq.\ (\ref{y_evol}). Assume that the stepsize
sequence satisfies $\alpha(k)\to 0$ as $k$ goes to infinity.
Then, for all $i$, we have
\[\lim_{k\to \infty} E[\|x_i(k)-y(k)\|] = 0,\quad \hbox{and}\]\[\liminf_{k\to \infty}
\|x_i(k)-y(k)\| = 0\qquad \hbox{with probability one.}\]
\label{mean-consensus-diffset}
\end{proposition}

\begin{proof}
From Lemma \ref{est-xkyk}, we have
\begin{eqnarray*}
\|x_i(k)-y(k)\| &\le & m\rho(k-1,0)\sum_{j=1}^m \|x_j(0)\| +
mL \sum_{r=0}^{k-2} \rho(k-1,r+1)\a(r) + 2\a(k-1) L\cr && +
\sum_{r=0}^{k-2} \rho(k-1,r+1)\sum_{j=1}^m \| e_j(r)\| +
\|e_i(k-1)\|+\frac{1}{m}\sum_{j=1}^m \|e_j(k-1)\|.
\end{eqnarray*}
Taking the expectation of both sides and using the estimate for
the disagreement metric $\rho(k,s)$ from Proposition
\ref{prop:compact-contraction}, i.e., for all $k\ge s\ge 0$,
\[E[\rho(k,s)]\le \kappa e^{-\mu(k-s)},\] for some scalars $\kappa,\mu>0$, we obtain
\begin{eqnarray*}
E[\|x_i(k)-y(k)\|] &\le & m \kappa e^{-\mu(k-1)}\sum_{j=1}^m
\|x_j(0)\| + mL\kappa \sum_{r=0}^{k-2}  e^{-\mu(k-r-2)}\a(r) +
2\a(k-1) L\cr && + \kappa\sum_{r=0}^{k-2}
e^{-\mu(k-r-2)}\sum_{j=1}^m \| e_j(r)\| +
\|e_i(k-1)\|+\frac{1}{m}\sum_{j=1}^m \|e_j(k-1)\|.
\end{eqnarray*}
By taking the limit superior in the preceding relation and
using the facts that $\alpha(k)\to 0$, and $\|e_i(k)\|\to 0$
for all $i$ as $k\to \infty$ (cf.\ Lemma \ref{error-behav}(b)),
we have for all $i$,
\begin{eqnarray*}
\limsup_{k\to \infty}E[\|x_i(k)-y(k)\|] \le  mL \kappa
\sum_{r=0}^{k-2} e^{-\mu(k-r-2)} \a(r) \ + \kappa \sum_{r=0}^{k-2}
e^{-\mu(k-r-2)}\sum_{j=1}^m \| e_j(r)\|.
\end{eqnarray*}
Finally, since $\sum_{k=0}^\infty e^{-\mu k}< \infty$ and both
$\alpha(k)\to 0$ and $\|e_i(k)\|\to 0$ for all $i$, by Lemma
\ref{lemma:seq}, we have
\[\lim_{k\to \infty}\sum_{r=0}^{k-2} e^{-\mu(k-r-2)}\a(r)=0\quad \hbox{and} \quad
\lim_{k\to \infty} \sum_{r=0}^{k-2} e^{-\mu(k-r-2)} \sum_{j=1}^m \|
e_j(r)\|=0.\] Combining the preceding two relations, we have
\[\lim_{k\to \infty}E[\|x_i(k)-y(k)\|]=0.\]
The second part of proposition follows using Fatou's Lemma and
a similar argument used in the proof of Proposition
\ref{mean-consensus-sameset}.
\end{proof}

The next proposition uses the compactness of the constraint
sets to strengthen this result and establish almost sure
consensus among the agent estimates.

\begin{proposition} Let Assumptions \ref{connectivity}, \ref{ass.stoch_weights},
\ref{ass.self_conf} and \ref{compactconst} hold.
 Let $\{x_i(k)\}$ be the sequence
generated by the algorithm
(\ref{convex-comb})-(\ref{proj-error}), and $\{y(k)\}$ be
defined in Eq.\ (\ref{y_evol}). Assume that the stepsize
sequence satisfies $\alpha(k)\to 0$. Then, for all $i$, we have
\[\lim_{k\to \infty} \|x_i(k)-y(k)\| = 0 \qquad \hbox{with
probability one}.\] \label{difconstasconsensus}
\end{proposition}

\begin{proof}
Using the iterations (\ref{subgradient-step}) and
(\ref{y_evol}), we obtain for all $k\ge 1$ and $i$,
\begin{eqnarray*}
y(k+1)-x_i(k+1) = \Big(y(k)-\sum_{j=1}^m a_{ij}(k) x_j(k)\Big)& -&
\alpha(k)\Big({1\over m} \sum_{j=1}^m d_j(k) -d_i(k)\Big)\\& +&
\Big({1\over m} \sum_{j=1}^m e_j(k) - e_i(k)\Big).
\end{eqnarray*}
Using the doubly stochasticity of the weights $a_{ij}(k)$ and
the subgradient boundedness (which holds by Assumption
\ref{compactconst}), this implies that
\begin{equation}\sum_{i=1}^m \|y(k+1)-x_i(k+1)\| \le \sum_{i=1}^m
\|y(k)-x_i(k)\| + 2Lm\alpha(k) + 2\sum_{i=1}^m
\|e_i(k)\|.\label{convseq}\end{equation}  Since $\a(k)\to 0$, it
follows from Lemma \ref{error-behav}(b) that $\|e_i(k)\|\to 0$
for all $i$. Eq.\ (\ref{convseq}) then yields
\begin{eqnarray*}
\limsup_{k\to \infty}\sum_{i=1}^m \|y(k+1)-x_i(k+1)\| &\le&
\liminf_{k\to \infty}\sum_{i=1}^m \|y(k)-x_i(k)\| \\ &&\ +
\lim_{k\to
\infty} \Big(2Lm\alpha(k) + 2\sum_{i=1}^m \|e_i(k)\|\Big)\\
&=& \liminf_{k\to \infty}\sum_{i=1}^m \|y(k)-x_i(k)\|.
\end{eqnarray*}
Using $x_i(k)\in X_i$ for all $i$ and $k$, it follows from
Assumption \ref{compactconst} that the sequence $\{x_i(k)\}$ is
bounded for all $i$. Therefore, the sequence $\{y(k)\}$
[defined by $y(k) = {1\over m}\sum_{i=1}^m x_i(k)$, see Eq.\
(\ref{yxaver})], and also the sequences $\|y(k)-x_i(k)\|$ are
bounded. Combined with the preceding relation, this implies
that the scalar sequence $\sum_{i=1}^m \|y(k)-x_i(k)\|$ is
convergent.

By Proposition \ref{mean-consensus-diffset}, we have
\[\liminf_{k\to \infty} \|x_i(k)-y(k)\| = 0\qquad \hbox{with probability one.}\]
Since $\sum_{i=1}^m \|y(k)-x_i(k)\|$ converges, this implies
that for all $i$,
\[\lim_{k\to \infty}\|x_i(k)-y(k)\|=0\qquad \hbox{with probability one},\]
completing the proof.
\end{proof}

The next theorem states our main convergence result for agent
estimates when the constraint sets are different under some
assumptions on the stepsize rule.

\begin{theorem} Let Assumptions \ref{connectivity}, \ref{ass.stoch_weights},
\ref{ass.self_conf}  and \ref{compactconst} hold. Let
$\{x_i(k)\}$ be the sequence generated by the algorithm
(\ref{convex-comb})-(\ref{proj-error}). Assume that the
stepsize sequence satisfies $\sum_k\alpha(k)=\infty$ and
$\sum_k \alpha^2(k)<\infty$. Then, there exists an optimal
solution $x^*\in X^*$ such that for all $i$
\[\lim_{k\to \infty} x_i(k)  = x^*\qquad \hbox{with probability one}.\]
\end{theorem}

\begin{proof}
From Lemma \ref{key-rel}(b), we have for some $z^*\in X^*$,
\begin{equation}\sum_{i=1}^m\|x_i(k+1)-z^*\|^2 \le \sum_{i=1}^m\|x_i(k)-z^*\|^2
+\a^2(k) \sum_{i=1}^m \|d_i(k)\|^2 -2\a(k) \sum_{i=1}^m
(f_i(v_i(k))-f_i(z^*)).\label{window-rel}\end{equation} We show that
the preceding implies that \begin{equation}\liminf_{k\to
\infty} \sum_{i=1}^m f_i(v_i(k))\le
f(z^*)=f^*.\label{optvalue}\end{equation} Suppose to arrive at
a contradiction that $\liminf_{k\to \infty} \sum_{i=1}^m
f_i(v_i(k))>f^*$. This implies that there exist some $K$ and
$\epsilon>0$ such that for all $k\ge K$, we have
\[\sum_{i=1}^m f_i(v_i(k))> f^* +\epsilon.\]
Summing the relation (\ref{window-rel}) over a window from $K$
to $N$ with $N>K$, we obtain
\[\sum_{i=1}^m\|x_i(N+1)-z^*\|^2 \le
\sum_{i=1}^m\|x_i(K)-z^*\|^2 +m L^2\sum_{k=K}^N \a^2(k) -2\epsilon
\sum_{k=K}^N\a(k).\] Letting $N\to \infty$, and using
$\sum_k\alpha(k)=\infty$ and $\sum_k \alpha^2(k)<\infty$, this
yields a contradiction and establishes the relation in Eq.\
(\ref{optvalue}).

By Proposition \ref{difconstasconsensus}, we have
\begin{equation}\lim_{k\to \infty} \|x_i(k)-y(k)\|=0\qquad
\hbox{with probability one}.\label{est-cons}\end{equation} Since
$v_i(k)=\sum_{j=1}^m a_{ij}(k) x_j(k)$, using the stochasticity
of the weight vectors $a_i(k)$, this also implies
\begin{equation}\lim_{k\to \infty} \|v_i(k)-y(k)\|\le \lim_{k\to \infty}
\sum_{j=1}^m a_{ij}(k) \|x_j(k)-y(k)\|=0\qquad \hbox{with
probability one}.\label{vconv}\end{equation} Combining Eqs.\
(\ref{optvalue}) and (\ref{vconv}), we obtain
\begin{equation}\liminf_{k\to \infty} f(y(k)) \le f^*\qquad
\hbox{with probability one}.\label{limopt}\end{equation}

From Lemma \ref{error-behav}(a), the sequence
$\{\sum_{i=1}^m\|x_i(k)-z^*\|\}$ is convergent for all $z^*\in
X^*$. Combined with Eq.\ (\ref{est-cons}), this implies that
the sequence $\{\|y(k)-z^*\|\}$ is convergent with probability
one. Therefore, $y(k)$ is bounded and it must have a limit
point. Moreover, since $x_i(k) \in X_i$ for all $k\ge 0$ and
$X_i$ is a closed set, all limit points of the sequence
$\{x_i(k)\}$ must lie in the set $X_i$ for all $i$. In view of
Eq.\ (\ref{est-cons}), this implies that all limit points of
the sequence $\{y(k)\}$ belong to the set $X$. Hence, from Eq.\
(\ref{limopt}), we have
\[\liminf_{k\to \infty} f(y(k)) = f^*\qquad
\hbox{with probability one}.\] Using the continuity of $f$ (due to
convexity of $f$ over $\mathbb{R}^n$), this implies that one of
the limit points
 of $\{y(k)\}$ must belong to $X^*$; denote this limit
point by $x^*$. Since the sequence $\{\|y(k)-x^*\|\}$ is
convergent, it follows that $y(k)$ can have a unique limit
point, i.e., $\lim_{k\to \infty} y(k)=x^*$ with probability
one. This and $\lim_{k\to\infty}\|x_i(k)-y(k)\|=0$ with
probability one imply that each of the sequences $\{x_i(k)\}$
converges to the same $x^*\in X^*$ with probability one.
\end{proof}

\section{Conclusions}\label{sec:conclusions}

We studied distributed algorithms for multi-agent optimization
problems over randomly-varying network topologies. We adopted a
state-dependent communication model, in which the availability of
links in the network is probabilistic with the probability dependent
on the agent states. This is a good model for a variety of
applications in which the state represents the position of the
agents (in sensing and communication settings), or the beliefs of
the agents (in social settings) and the distance of the agent states
affects the communication and information exchange among the agents.

We studied a projected multi-agent subgradient algorithm for this
problem and presented a convergence analysis for the agent
estimates. The first step of our analysis establishes convergence
rate bounds for a disagreement metric among the agent estimates.
This bound is time-nonhomogeneous in that it depends on the initial
time. Despite this, under the assumption that the stepsize sequence
decreases sufficiently fast, we proved that agent estimates converge
to an almost sure consensus and also to an optimal point of the
global optimization problem under some assumptions on the constraint
sets.

The framework introduced in this paper suggests a number of
interesting further research directions. One future direction is to
extend the constrained optimization problem to include both local
and global constraints. This can be done using primal algorithms
that involve projections, or using primal-dual algorithms where dual
variables are used to ensure feasibility with respect to global
constraints. Another interesting direction is to consider different
probabilistic models for state-dependent communication. Our current
model assumes the probability of communication is a continuous
function of the $l_2$ norm of agent states. Considering other norms
and discontinuous functions of agent states is an important
extension which is relevant in a number of engineering and social
settings.

\newpage

\bibliographystyle{amsplain}
\bibliography{distributed}

\end{document}